\documentclass[12pt]{article}

\usepackage{graphicx}
\usepackage{amssymb}
\usepackage{amsmath}
\usepackage{color}

\newtheorem{theorem}{Theorem}
\newtheorem{statement}{Statement}
\newtheorem{lemma}{Lemma}

\newtheorem{corollary}{Corollary}
\newtheorem{remark}{Remark}

\newcommand{\kom}[1]{}
\renewcommand{\kom}[1]{{\bf [#1]}}
\definecolor{blau}{rgb}{0.1,0.0,0.9}

\newcounter{komcounter}
\numberwithin{komcounter}{section}

\begin{document}

\title{\bf Estimates of size of cycle in a predator-prey system}
\author{Niklas L.P. Lundstr\"{o}m$^{\tiny 1}$ and
Gunnar S\"oderbacka$^{\tiny 2}$
 \linebreak \\\\
$^{\tiny 1}$\it \small Department of Mathematics and Mathematical Statistics, Ume{\aa} University\\
\it \small SE-90187 Ume{\aa}, Sweden\/{\rm ;}
\it \small niklas.lundstrom@umu.se\linebreak \\\\
$^{\tiny 2}$\it \small {\AA}bo Akademi, 20500 {\AA}bo, Finland\/{\rm ;}
\it \small gsoderba@abo.fi \linebreak \\\\}

\maketitle

\begin{abstract}
We consider a Rosenzweig-MacArthur predator-prey system which incorporates logistic growth of the prey
in the absence of predators and a Holling type II functional response for interaction between predators and preys.
We assume that parameters take values in a range which guarantees that all solutions tend to a unique limit cycle and prove estimates for the maximal and minimal predator and prey population densities of this cycle.
Our estimates are simple functions of the model parameters and
hold for cases when the cycle exhibits small predator and prey abundances and large amplitudes.
The proof consists of constructions of several Lyapunov-type functions and derivation of a large number of non-trivial estimates which are also of independent interest.\\


\end{abstract}


\setcounter{equation}{0} \setcounter{theorem}{0}

\section{Introduction}

The dynamical relationship between predators and preys,
most simply described by Lotka-Volterra-type ordinary differential equations,
has been investigated widely in recent years.
One well known mathematical model describing this relationship is the
Rosenzweig-MacArthur extension of the classical Lotka-Volterra model, see e.g. \cite{R71,T03,M74,Mu89,Y89,H98},
in which various interaction rates between the populations have nonlinear dependence on the prey concentration according to
\begin{align}\label{Rosenzweig-MacArthur}
 \frac{dS}{dt} &\,=\, r S \left(1-\frac{S}{K}\right) -  \frac{q X S}{H + S}, \notag\\
 \frac{dX}{dt} &\,=\,  \frac{p X S}{H + S} - d X.
\end{align}
Here, $S \,=\, S\left(t\right)$ and $X \,=\, X\left(t\right)$ denotes the population densities of prey and predator,
respectively, and $r, K, q, H, p$ and $d$ are positive parameters.
The biological meanings of the parameters are the following:
$r$ is the intrinsic growth rate of the prey;
$K$ is the prey carrying capacity;
$q$ is the maximal consumption rate of predators;
$H$ is the amount of prey needed to achieve one-half of $q$;
$d$ is the per capita death rate of predators;
and $p$ is the efficiency with which predators convert consumed prey into new predators.

In this paper, we prove analytical estimates of the size of a limit cycle in the following version of system \eqref{Rosenzweig-MacArthur}:
\begin{align}\label{pp}
\frac{ds}{d\tau} &\,=\, \left(h\left(s\right)-x\right)\,s,  \notag\\
\frac{dx}{d\tau} &\,=\, \left(s-\lambda\right)\,x, \qquad \text{where}\qquad h\left(s\right)\,=\,\left(1-s\right)\left(s+a\right)
\end{align}
and 
$s \,=\, s\left(\tau\right)$ and $x \,=\, x\left(\tau\right)$ denote the population densities of prey and predator,
respectively.
We will focus on the dynamics of system \eqref{pp} when the parameters $a$ and $\lambda$ take on small values,
namely, we assume
\begin{equation}\label{pk}
 a \,<\, 0.1 \qquad \mbox{and} \qquad \lambda \,<\, 0.1.
\end{equation}

In order to describe the simple relation between the above Rosenzweig-MacArthur system in \eqref{pp}
and the more standard version given in \eqref{Rosenzweig-MacArthur},
we observe that by introducing
the scaled time $\tau$,
the state variables $s \,=\, s\left(\tau\right)$ and $x \,=\, x\left(\tau\right)$
and the parameters $a$, $b$ and $\lambda$ according to
\begin{align*}
\tau \,=\, \int \frac{r K}{H + S\left(t\right)} dt, \quad s \,=\, \frac{S}{K}, \quad x \,=\, \frac{q X}{r K}, \quad a \,=\, \frac{H}{K},  \quad b \,=\, \frac{p - d}{r} \quad \text{and} \quad \lambda \,=\, \frac{d H}{r K},
\end{align*}
the standard system in \eqref{Rosenzweig-MacArthur} transforms to system \eqref{pp} when $b \,=\, 1$.

Rosenzweig-MacArthur systems incorporate logistic growth of the prey
in the absence of predators and a Holling type II functional response
(Michaelis-Menten kinetics) for interaction between predators and preys.
A literature survey shows that the model has been widely used in real life ecological applications,
see e.g. \cite{G-OR-J03,WWS16,M72,RM63,R71},
including the spatiotemporal dynamics of an aquatic community of phytoplankton and zooplankton \cite{MPTML02}
as well as dynamics of microbial competition \cite{HHW77,BHW83,SW95}.

From a mathematical point of view, the dynamics of systems of type \eqref{Rosenzweig-MacArthur} and \eqref{pp} has been frequently studied, see e.g. \cite{H77,C81,K83,D88,H88,KF88,HM89,HH95,EOS96} and the references therein.
In particular, system \eqref{pp} always has a unique positive equilibrium at
$\left(x,s\right) \,=\, \left(\left(1-\lambda\right)\left(\lambda + a\right), \lambda\right)$
which attracts the whole positive space when $2\lambda + a \,>\, 1$.
At $2\lambda + a \,=\, 1$ there is a Hopf bifurcation in which the equilibrium loses stability and a stable limit cycle, surrounding the equilibrium, is created.
In particular, for $2\lambda + a \,<\, 1$ the equilibrium is a source and the cycle attracts the whole positive space (except the source), \cite{C81}.

Our main results are analytical estimates of the size of this unique limit cycle when
parameters values of $a$ and $\lambda$ are small.
Namely, we assume \eqref{pk} and in such cases the amplitude of the cycle becomes large and $x$ and $s$ may become very small during a large portion of the cycle.
Biologically, this means that the modeled population exhibits small predator and prey abundances,
indicating that the population suffers a relatively high risk of going extinct because of random perturbations.
This underscores the importance of understanding the dynamics of systems of type \eqref{pp} under assumption \eqref{pk}.
To further motivate our analytical estimates, we mention that it is notrivial to obtain accurate numerical results by integrating the equations \eqref{pp} using standard numerical methods when $a$ and $\lambda$ are small,
see Section \ref{sec:numerics}.

Before stating our main results,
let us note that the above Rosenzweig-MacArthur systems are very simplified models of reality and therefore
usually not directly applicable in biology without modifications.
For example, it is clear from our main results that in model \eqref{pp} predator and prey populations can decrease to unacceptable low abundances and still survive.
However, we believe that even though our main results are proved for such simple models,
they may be useful when investigating dynamics also in more complex and realistic systems,
such as, e.g., systems modeling the interactions of several predators and one prey,
or seasonally dependent systems,
see e.g. \cite{BHW83,RMK93,EOS96}.

Our main results are summarized in the following theorem.

\begin{theorem}\label{th:main}
Let $x_{max}$ and $s_{max}$ be the maximal $x$- and $s$-values and let $x_{min}$ and $s_{min}$ be
the minimal $x$- and $s$-values in the unique limit cycle of system \eqref{pp} under assumption (\ref{pk}).
Then the predator density satisfies
\begin{align*}
1 \,<\, x_{max} \,<\, 1.6 \qquad \text{and} \qquad \exp{\left(-\frac{x_{max}}{a}\right)} \,<\,  x_{min} \,<\, \exp{\left(-\frac{x_{max}}{a\kappa_1}\right)},
\end{align*}
and the prey density satisfies
\begin{align*}
0.9 \,<\, s_{max} \,<\, 1 \qquad \text{and} \qquad \exp{\left(-\frac{ x_{max}}{\lambda\kappa_2}\right)} \,<\, s_{min} \,<\, \exp{\left(-\frac{ x_{max}}{\lambda\kappa_3}\right)},
\end{align*}
where
\begin{align*}
1
\,<\, \kappa_1
\,=\, \frac{1 + e^{-2}\frac{\lambda}{a}}{1 - 2.1 \lambda - 0.31 a}
\,<\, 1.32 \left(1 + e^{-2}\frac{\lambda}{a}\right),
\end{align*}
%
%
\begin{equation*}
1 \,>\, \kappa_2 \,=\, \frac{1}{1 +\lambda \left(1-\ln \left(\lambda\right)\right)} \,>\, 0.75 \quad \textrm{and} \quad 1 \,<\, \kappa_3 \,=\, \frac{1}{1 - 0.3\lambda - a \left(1.3 - \ln \left(a\right)\right)} \,<\, 1.64.
\end{equation*}
\end{theorem}
From the expressions for $\kappa_1, \kappa_2$ and $\kappa_3$ in Theorem \ref{th:main} we conclude that
\begin{align*}
&\kappa_1 \searrow 1 \qquad \textrm{if} \qquad a \to 0 \quad \textrm{and} \quad \frac{\lambda}{a} \to 0,\\
&\kappa_2 \nearrow 1 \qquad \textrm{if} \qquad \lambda \to 0,\\
&\kappa_3 \searrow 1 \qquad \textrm{if} \qquad a \to 0 \quad \textrm{and} \quad \lambda \to 0.
\end{align*}
Therefore, Theorem \ref{th:main} yields the following remark.

\begin{remark}\label{re:approx_main}
Let $x_{max}$ and $s_{max}$ be the maximal $x$- and $s$-values and let $x_{min}$ and $s_{min}$ be
the minimal $x$- and $s$-values in the unique limit cycle of system \eqref{pp} under assumption \eqref{pk}.
If both $a$ and $\frac{\lambda}{a}$ are small, then the estimate
$$
x_{min} \,\approx\, \exp\left({-\frac{x_{max}}{a}}\right)
$$
is good for the minimal predator biomass of the unique limit cycle.
Similarly, if both $a$ and $\lambda$ are small, then the estimate
$$
s_{min} \,\approx\, \exp\left({-\frac{x_{max}}{\lambda}}\right)
$$
is good for the minimal prey biomass of the unique limit cycle.
\end{remark}

Before discussing the outline of the proof of Theorem \ref{th:main},
we state its analogue for the more standard version of the Rosenzweig-MacArthur system given in \eqref{Rosenzweig-MacArthur} as a corollary.
In this setting, assumption \eqref{pk} takes the form
\begin{align} \label{Main assumptions*}
\frac{H}{K} \,<\, 0.1, \qquad \frac{d H}{r K} \,<\, 0.1  \qquad \text{and} \qquad \frac{p - d}{r} \,=\, 1.
\end{align}
The biological meaning of the first two inequalities is that the half saturation rate for predators ($H$) is assumed small compared to the carrying capacity of the prey ($K$),
and that the death rate of predators ($d$) is assumed small compared to the growth rate of the prey ($r$) times $K/H$.
The third assumption in \eqref{Main assumptions*} says that the growth rate of the prey ($r$) equals the difference between
the efficiency of the predators ($p$) and the death rate of predators ($d$).
Theorem \ref{th:main} immediately implies the following result.

\begin{corollary}\label{th:main*}
Let $X_{max}$ and $S_{max}$ be the maximal predator and prey densities
and let $X_{min}$ and $S_{min}$ be the minimal predator and prey densities
of the unique limit cycle in system \eqref{Rosenzweig-MacArthur} under assumption \eqref{Main assumptions*}.
Then the predator density satisfies
%
%
\begin{align*}
1 \,<\, \frac{q X_{max}}{r K} \,<\, 1.6 \qquad \text{and} \qquad
\exp{\left(-\frac{q X_{max}}{r H}\right)} \,<\, \frac{q X_{min}}{r K} \,<\,  \exp{\left(-\frac{q X_{max}}{r H \kappa_1}\right)},
\end{align*}
and the prey density satisfies
\begin{align*}
0.9 \,<\, \frac{S_{max}}{K} \,<\, 1 \qquad \text{and} \qquad
\exp{\left(-\frac{q X_{max}}{d H \kappa_2}\right)} \,<\, \frac{S_{min}}{K} \,<\, \exp{\left(-\frac{q X_{max}}{d H\kappa_3}\right)},
\end{align*}
where $\kappa_1$, $\kappa_2$ and $\kappa_3$ are given by Theorem \ref{th:main} with
$a \,=\, \frac{H}{K}$ and $\lambda \,=\, \frac{d H}{r K}$.

Moreover, if $\frac{H}{K}$ and $\frac{d}{r}$ are small,
then the estimate $X_{min} \approx \frac{r K }{q}\exp{\left(-\frac{q X_{max}}{r H}\right)}$
is good for the minimal predator biomass,
and if $\frac{H}{K}$ and $\frac{dH}{rK}$ are small,
then the estimate $S_{min} \approx K \exp{\left(-\frac{ q X_{max}}{d H }\right)}$ is good for the minimal prey biomass.
\end{corollary}

The proof of Theorem \ref{th:main} consists of constructions of several Lyapunov-type functions and derivation of a large number of non-trivial estimates.
We believe that these methods and constructions have values also beyond this paper
as they present methods and ideas that, potentially, can be useful for proving analogous results
for dynamics in similar systems as well as in more complex systems.

The proof is constructed in a way such that Theorem \ref{th:main} is a direct consequence of four statements,
namely Statement 1-4,
which we prove in Sections \ref{sec:Region1}, \ref{sec:Region23} and \ref{sec:Region4}.
In addition to the estimates in Theorem \ref{th:main}
it is also possible to find, from these statements, a positively invariant region trapping the unique limit cycle inside.
In fact, the limit cycle will be inside an outer boundary consisting of the part of a trajectory $\hat T$ with initial condition $x\left(0\right)\,=\,1.6, \, s\,=\,\lambda$ and the part of $s\,=\,\lambda$ between
$\left(1.6, \lambda\right)$ and the next intersection with $s\,=\,\lambda$ when $x\,>\,h\left(\lambda\right)$.
It will also be outside an inner  boundary consisting of the part of a trajectory $\check T$ with initial condition $x\left(0\right)\,=\,1, \, s\,=\,\lambda$ and the part of $s\,=\,\lambda$ between
$\left(1, \lambda\right)$ and the next intersection with $s\,=\,\lambda$ when $x\,>\,h\left(\lambda\right)$.
Estimates for these boundaries can be found from given statements and lemmas,
even if we do not write them explicitly here.
We also point out that better
but more complicated
estimates than those summarized in Theorem \ref{th:main} follow from lemmas which are used for the
proof of Statements 1-4 and Theorem 1.

To outline the proof of Theorem \ref{th:main} we first observe that the coordinate axes are invariant, and hence the region $x,s \,>\, 0$ is also invariant.
Therefore, we consider solutions only for positive $s$ and $x$.
Moreover, system \eqref{pp} has isoclines at $x \,=\, h\left(s\right)$ and $s \,=\, \lambda$,
which leads us to split the proof by introducing the following four regions:
\begin{itemize}
\item[] {Region 1}, where $x\,>\,h\left(s\right),\ s\,>\,\lambda$ and $x$ is growing and $s$ decreasing.
\item[] {Region 2}, where $x\,>\,h\left(s\right),\ s\,<\,\lambda$ and both $x$ and $s$ decrease.
\item[] {Region 3}, where $x\,<\,h\left(s\right),\ s\,<\,\lambda$ and $x$ decreases and $s$ grows.
\item[] {Region 4}, where $x\,<\,h\left(s\right),\ s\,>\,\lambda$ and both $x$ and $s$ increase.
\end{itemize}
\begin{figure}[h]
\begin{center}
\includegraphics[scale = 0.6]{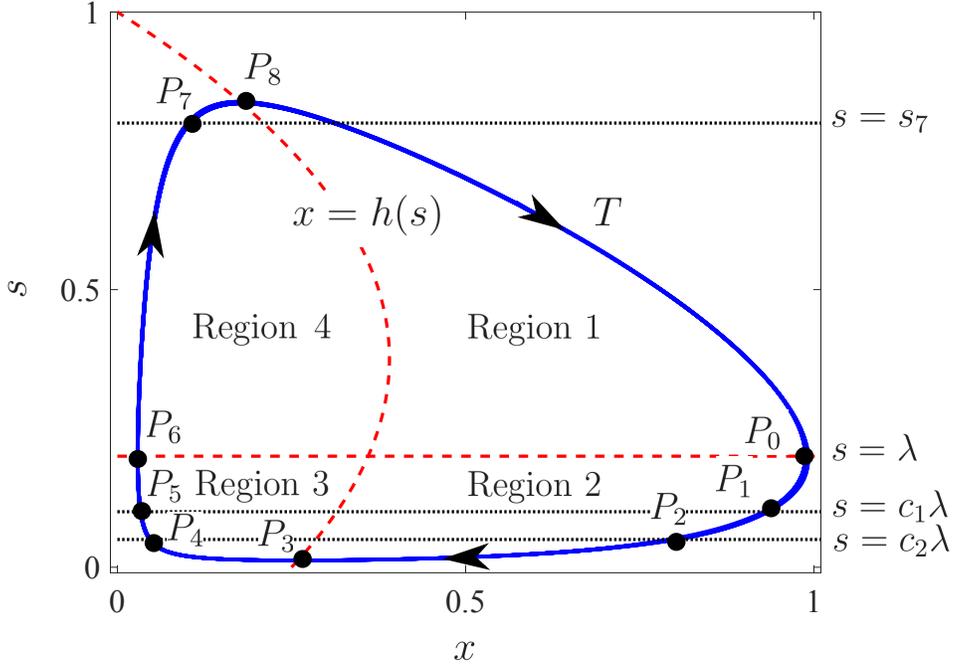}
\caption{Notations of the four Regions 1-4, the points $P_1$-$P_8$ on a trajectory $T$ (blue, solid),
and the isoclines $x \,=\, h(s)$ and $s \,=\, \lambda$ (red, dotted) of system \eqref{pp}.} 
\label{figdok}
\end{center}
\end{figure}
Any trajectory starting in Region 1 will enter Region 2 from where it will enter Region 3 and then Region 4 and finally Region 1 again,
and the behaviour repeats infinitely.
Figure \ref{figdok} illustrates the four regions together with isoclines and points which will be used in the proof of Theorem \ref{th:main}.
Behaviour and estimates for trajectories in different regions are examined in different sections.
Behaviour in Region 1 are examined in Section \ref{sec:Region1}.
Section \ref{sec:Region23} considers Regions 2 and 3 while Region 4 is considered in Section \ref{sec:Region4}.
The main results in Regions 1-4 will be concluded in Statements 1-4.
We end the paper by giving some numerical results in Section~\ref{sec:numerics}.



\setcounter{equation}{0} \setcounter{theorem}{0}

\section{Estimates in Region 1}
\label{sec:Region1}

We begin this section by proving a lemma which gives a bounded region
into which all trajectories will enter after sufficient time
and which will be used in several places in the proof of Theorem \ref{th:main}.

\begin{lemma}\label{le:lemma1}
Consider the function
$$
V_g\left(s\right) \,=\, \frac{\alpha \left(1-s\right)}{1+\beta \left(1-s\right)}
\qquad \text{where} \qquad
\alpha \,=\, 2-\lambda +a
\qquad \text{and} \qquad
\beta \,=\,\frac{\lambda +1}{a-2\lambda +3}.
$$
All solutions of system \eqref{pp} under condition \eqref{pk}
with positive initial values will enter into the region determined by the inequalities
$x \,<\, V_g\left(s\right)$, $x \,>\, 0$ ,$s \,>\, 0$ and remain there.
\end{lemma}




\noindent
{\it Proof of Lemma \ref{le:lemma1}}.
\begin{figure}[h]
\begin{center}
\includegraphics[scale = 0.5]{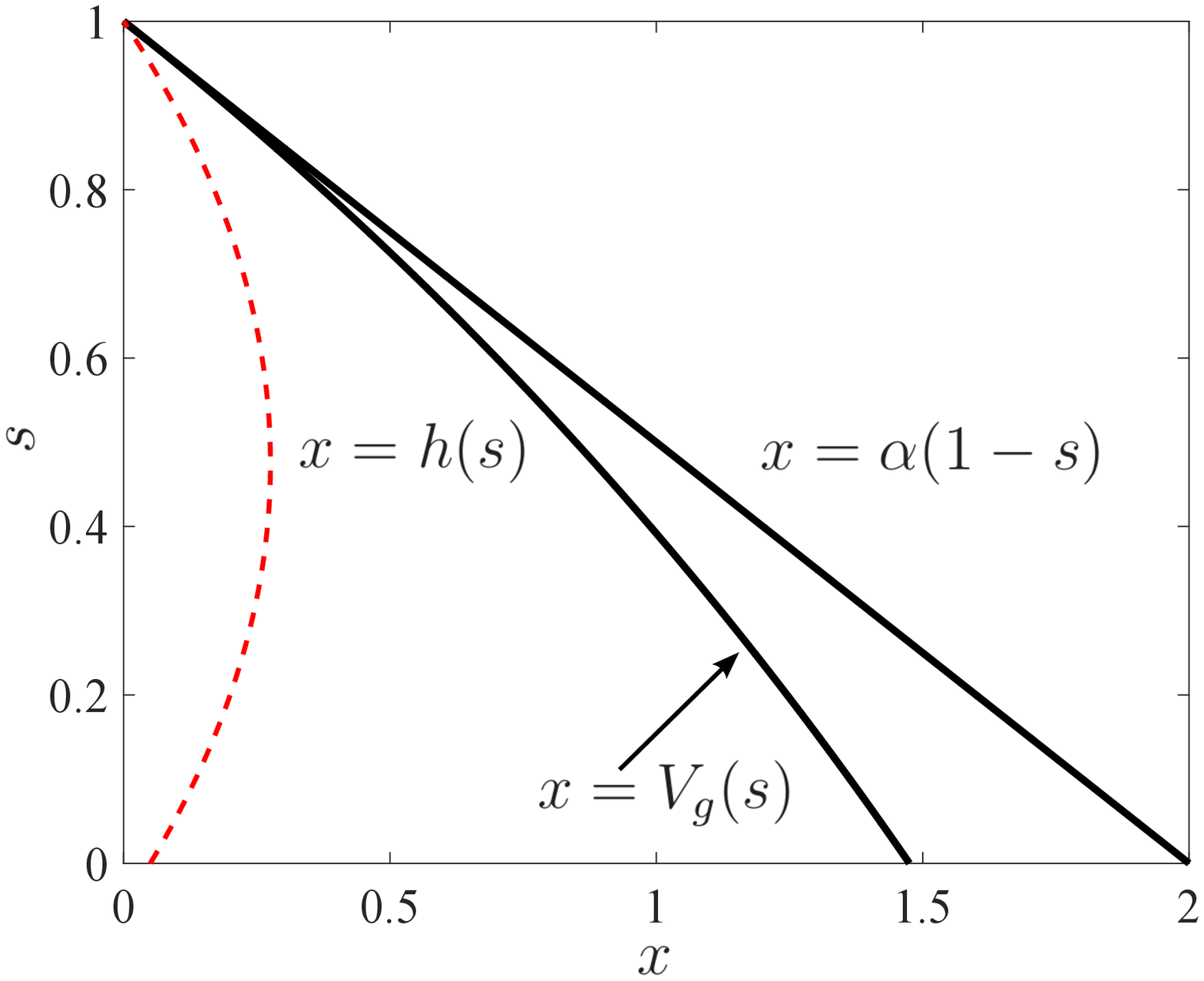}
\caption{Geometry in the proof of Lemma \ref{le:lemma1}. }
\label{fig:Lemma1}
\end{center}
\end{figure}
Let $U\,=\,x-\alpha \left(1-s\right)$.
Differentiation with respect to time and using \eqref{pp} yield
$$
U'\,=\,-\left(1-s\right)^2 \left(s+\lambda \right)\,\alpha +\gamma\, \left(\lambda s -s -as -\lambda\right)\,<\,0
$$
\noindent
for $\gamma \,\geq\, 0$ when $x\,=\,\alpha \left(1-s\right) +\gamma $.
Thus all trajectories will enter the region $x \,<\, \alpha \left(1-s\right)$ and remain there, see Figure \ref{fig:Lemma1}.
Let also
$$
V\,=\,\left(1+\beta \left(1-s\right)\right)x -\alpha \left(1-s\right),
$$
and notice that $\frac{\partial V}{\partial x} , \frac{\partial V}{\partial s} \,>\, 0$ since $0 \,<\, \beta \,<\, 1 \,<\, \alpha$ and, after sufficient time, $s \,<\, 1$ and $x \,<\, \alpha$.
Calculating the derivative of $V$ with respect to time we get
\[V'\,=\,-\frac{\alpha  \,{\left( 1-s\right) }^{3}\,S }{{\left( 4-\lambda\,s-s+a-\lambda\right) }^{2}}\]
\noindent
at $x \,=\, V_g$, where
\[
S\,=\,\left( \lambda+1\right) \,\left( a-3\,\lambda+2\right) \,s+\lambda\,{\left( a-\lambda+4\right) }^{2}.
\]
Because $S\,>\,0$ for $a,\lambda \,<\,0.1$ we get $V'\,\leq\, 0$ for $0\,<\,s\,\leq\, 1$ and $V'\,=\,0$ only for $s\,=\,1$.
Because $\beta\,>\,0$ we have $V_g\,<\,\alpha \left(1-s\right)$ and since $h\left(s\right) \,<\, V_g\left(s\right)$ for $0\,<\,s\,<\,1$ all trajectories entering
$x\,<\,\alpha \left(1-s\right)$ also enter region $x\,<\,V_g$, where they remain because of the sign of $V'$. $\hfill\Box$\\


The maximal $x$-value for a trajectory is attended when it escapes from Region~1 to Region 2. In this section we will give estimates for maximal $x$-value, when trajectory starts on boundary of Region 1.

\begin{statement}
Any trajectory starting on the isocline $x \,=\, h\left(s\right), s \,>\, \lambda$ has a maximum $x_0$ before it enters Region 2 and $x_0 \,<\, 1.6$. Moreover, if the trajectory starts from a point where $s \,>\, 0.9$, then $x_0 \,>\, 1$.
\end{statement}

\noindent
We formulate the last part of the statement as lemma with own proof.

\begin{lemma}\label{le:lemma2}
Any trajectory starting on the isocline $x\,=\,h\left(s\right), s \,>\, 0.9$ has a maximum $x_0$ before it enters Region 2 and $x_0\,>\,1$.
\end{lemma}

\noindent
{\it Proof of Lemma \ref{le:lemma2}}.
In Region 1 the $x$-value on the trajectory is growing while the $s$-value is decreasing,
and $x'$ is smallest for greatest $\lambda$ and $s'$ is smallest for smallest $a$.
This implies that in Region 1, for any $a, \lambda \in [0,0.1)$,
the $x$-value for a trajectory of system \eqref{pp} is always growing stronger than the $x$-value for a trajectory of the system obtained for $a\,=\,0$ and $\lambda \,=\, 0.1$, since $\vert\frac{dx}{ds}\vert$ will then be smallest.
By this fact we are able to construct a bound for the minimal value of $x_0$ by using system \eqref{pp} with $a\,=\,0$ and $\lambda \,=\, 0.1$ fixed.

We define the continuous function $f$ by

\[ f\left(s\right)\,=\,\left\{ \begin{array}{ll}
           0.513+1.33 s - 2 s^2 & 0.7\,<\,s\,\leq\, 0.9 \\
           1.045-0.13s-s^2 & 0.5\,<\,s\,\leq\, 0.7 \\
           1.08-0.2s-s^2  & 0.3\,<\,s\,\leq\, 0.5 \\
           0.975+0.45s-2s^2 & 0.1\,\leq\, s\,\leq\, 0.3
            \end{array}
          \right.
\]

\noindent
and consider the function $Y$ defined by $Y\left(x,s\right)\,=\,x-f\left(s\right)$.
The derivative of $Y$ with respect to time,
considering system \eqref{pp} with $a\,=\,0,\ \lambda\,=\,0.1$ and substituting $x\,=\,f\left(s\right)$,
is a fourth order polynomial on each piece of definition.
By standard techniques it can be shown that this derivative is positive on each piece.
Thus, $x$ grows faster than $f\left(s\right)$ on the curve $x\,=\,f\left(s\right)$.
\begin{figure}[h]
\begin{center}
\includegraphics[scale = 0.5]{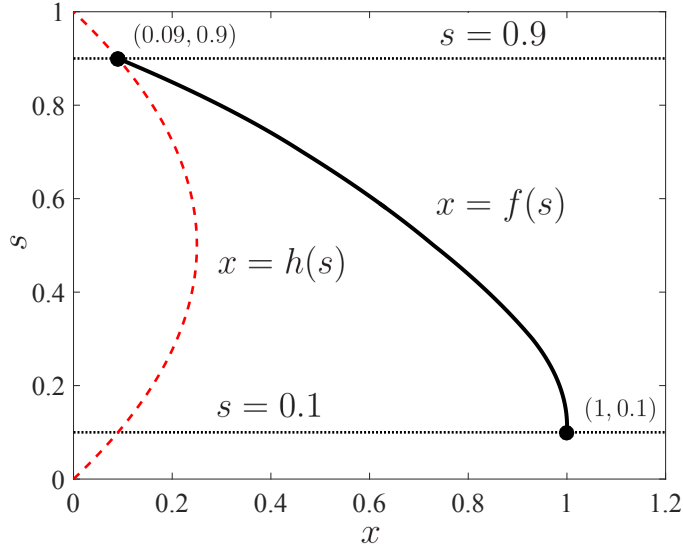}
\caption{Geometry in the proof of Lemma \ref{le:lemma2}. }
\label{fig:Lemma2}
\end{center}
\end{figure}
Moreover, for $s\,=\,0.9$ we get $f\left(s\right)\,=\,0.09$,
meaning that the point $\left(f\left(s\right),s\right) \,=\, \left(0.09, 0.9\right)$ is on the isocline $x \,=\, h\left(s\right)$ because $h\left(0.9\right) \,=\, 0.09$,
see Figure \ref{fig:Lemma2}.
We conclude that trajectories intersect the pieces of $x\,=\,f\left(s\right)$ transversally
going from the region defined by $x\,<\,f\left(s\right)$ to region where $x\,>\,f\left(s\right)$.
The isocline of any system \eqref{pp} under condition \eqref{pk} is above the isocline for the system we considered,
meaning $x$ is greater and also trajectories cannot intersect $x\,=\,f\left(s\right)$ before $s\,<\,0.1$.
Moreover $f\left(0.1\right)\,=\,1$.
Thus any trajectory for any $a,\lambda \in [0,0.1)$ under our conditions that start on the isocline $x\,=\,h\left(s\right)$, $s \,>\, 0.9$,
will at $s\,=\,0.1$ have an $x$-value greater than 1 and consequently this holds also at $s\,=\,\lambda$.
Therefore, $x_0 \,>\, 1$ and the proof of Lemma \ref{le:lemma2} is complete. $\hfill\Box$\\

\noindent
{\it Proof of Statement 1}.
We first recall the notations from Lemma \ref{le:lemma1} and also the fact that
the maximum of a trajectory taken for $s \,=\, \lambda$ before it enters Region 2 is less than $V_g\left(\lambda\right)$,
and

\[
V_g\left(\lambda\right)\,=\,-\frac{\left( \lambda-1\right) \,\left( a-2\,\lambda+3\right) \,\alpha }{a-{\lambda}^{2}-2\,\lambda+4}.
\]

\noindent
For $a,\lambda \in (0, 0.1]$,
the derivative of $V_g\left(\lambda\right)$ with respect to $a$ is

\[
\frac{\left( 1-\lambda\right) \,\left( {a}^{2}-2\,{\lambda}^{2}\,a-4\,\lambda\,a+8\,a+3\,{\lambda}^{3}-{\lambda}^{2}-15\,\lambda+14\right) }{{\left( a-{\lambda}^{2}-2\,\lambda+4\right) }^{2}}\,>\,0,
\]

\noindent
and derivative of $V_g\left(\lambda\right)$ with respect to $\lambda$ is

\[
-\frac{{a}^{3}+{\lambda}^{2}\,{a}^{2}-8\,\lambda\,{a}^{2}+10\,{a}^{2}+20\,{\lambda}^{2}\,a-52\,\lambda\,a+35\,a-2\,{\lambda}^{4}-8\,{\lambda}^{3}+55\,{\lambda}^{2}-84\,\lambda+40}{{\left( a-{\lambda}^{2}-2\,\lambda+4\right) }^{2}}\,<\,0.
\]

\noindent
Thus, $V_g\left(\lambda\right)$ is less than its value for $a \,=\, 0.1$ and $\lambda \,=\, 0$,
which is less than 1.588.
Since the estimate from below follows from Lemma \ref{le:lemma2}, the proof of Statement 1 is complete. $\hfill\Box$


\setcounter{equation}{0} \setcounter{theorem}{0}

\section{Estimates in Region 2 and Region 3}
\label{sec:Region23}

We consider a trajectory $T$ of system \eqref{pp} under condition \eqref{pk}
with initial condition $x\left(0\right)\,=\,x_0, \, s\left(0\right)\,=\,\lambda$, where $1\,<\,x_0\,<\,1.6$.
We suppose $c_1$ and $c_2$ are such that $0\,<\,c_2\,<\,c_1\,<\,\lambda$.
If $T$ intersects $s\,=\,c_2\lambda$ before escaping Region 2 we denote the point
of first intersection with $s\,=\,c_1\lambda$ by $P_1\,=\,\left(x_1, c_1\lambda \right)$ and the point of first intersection with   $s\,=\,c_2\lambda$ by $P_2\,=\,\left(x_2, c_2\lambda \right)$.
We denote the next intersection with the isocline $x\,=\,h\left(s\right)$ by $P_3\,=\,\left(x_3, s_3\right)$, where $x_3\,=\,h\left(s_3\right)$. The second intersection with $s\,=\,c_2\lambda$ we denote by $P_4\,=\,\left(x_4, c_2\lambda \right)$ and the second intersection with $s\,=\,c_1\lambda$  by $P_5\,=\,\left(x_5, c_1\lambda \right)$. The next intersection with $s\,=\,\lambda$ we denote by $P_6\,=\,\left(x_6, \lambda \right)$.
The lowest $s$-value  of the trajectory before it escapes to Region 4 will be at $P_3$
and the lowest $x$-value at $P_6$.
The notations are illustrated in Figure \ref{figdok}, where there are added also points used in Section 4.
We point out that trajectory $T$ is normally not a cycle, even though such case is illustrated in Figure \ref{figdok}.

The main results in this section are given in
Statements 2 and 3 which give main estimates in Regions 2 and 3.
Statement 2 gives a lower and upper bound for minimal $x$-value and Statement 3 gives  a lower and upper bound for minimal $s$-value of the part of the trajectory in Regions 2 and 3. Lemma \ref{le:lemma3} gives a better upper estimate for lowest $x$-value  which is needed also in Section 4.  These estimates will also serve as upper and lower estimates for the unique cycle of system \eqref{pp} under condition \eqref{pk}.
In the proofs of Statement 2 and Lemma \ref{le:lemma3} we assume $c_1\,=\,e^ {-2},\, c_2\,=\,e^{-4}$.

We here give these three main results of this section.

\begin{lemma}\label{le:lemma3}
Trajectory $T$ intersects $s\,=\,e^{-4}\lambda$ before escaping Region 2 and for $x_6$ we have the estimate

\begin{equation}
x_6 \,<\, 1.015\, e^{-2A-\frac{\theta \left(\tilde x_2\right)}{a+e^{-4}\lambda}}
\label{ec3}
\end{equation}

\noindent
where

\begin{equation*}
\tilde x_2\,=\,x_0 -3.8\, \lambda, \qquad  \theta \left(\tilde x_2\right) \,=\, \tilde x_2 - \left(a+e^{-4}\lambda\right) \ln \left(\tilde x_2\right),\qquad 
A\,=\,\frac{\left(1-e^{-2}\right)\lambda}{a+e^{-2}\lambda}.
\end{equation*}
\end{lemma}

\noindent
More general estimates than in Lemma \ref{le:lemma3} and Statement 2 are given in Lemma \ref{le:lemma4} and \ref{le:lemma5}.
These are formulated for general choices of parameters $c_1$ and $c_2$, which are fixed in proofs of Lemma \ref{le:lemma3} and Statement 2.

\begin{statement}
For the intersection of trajectory $T$ with the isocline $s \,=\, \lambda$ at $P_6 \,=\, \left(x_6, s_6\right)$
the following estimates are valid for the $x$-value.
%
\begin{equation}
e^{-\frac{x_0}{a}} \,<\, x_6 \,<\, e^{-\frac{x_0}{a\kappa_1}},
\label{e3c}
\end{equation}
where
\begin{equation}
1
\,<\, \kappa_1
\,<\, \frac{1 + e^{-2}\frac{\lambda}{a}}{1 - 2.1 \lambda - 0.31 a}
\,<\, 1.32 \left(1 + e^{-2}\frac{\lambda}{a}\right).
\label{e3ck}
\end{equation}
\end{statement}

\noindent
From Statement 2 it follows that for small $\frac{\lambda}{a}$
and $a$ the estimate $e^{-\frac{x0}{a}}$ is good for the minimal $x$-value on trajectory $T$.

\begin{statement}
For the intersection of trajectory $T$ with the isocline $x \,=\, h\left(s\right)$ at $P_3 \,=\, \left(x_3, s_3\right)$ the following estimates are valid for the $s$-value. 
%
\begin{equation}
e^{-\frac{x_0}{\lambda\kappa_2}} \,<\,s_3 \,<\, e^{-\frac{x_0}{\lambda\kappa_3}},
\label{e3d}
\end{equation}
where
\begin{equation*}
1 \,>\, \kappa_2 \,>\, \frac{1}{1 +\lambda \left(1-\ln \left(\lambda\right)\right)} \,>\, 0.75,
\end{equation*}
and
\begin{equation*}
1 \,<\, \kappa_3 \,<\, \frac{1}{1 - 0.3\lambda - a \left(1.3 - \ln \left(a\right)\right)} \,<\, 1.64.
\end{equation*}
\end{statement}

\noindent
From Statement 3 we see that for small $\lambda,\, a$ the estimate $e^{-\frac{x0}{\lambda}}$ is good for the minimal $s$-value on $T$.

The proof of Statement 2 is following from Lemma \ref{le:lemma3} and Statement 1 and a short Lemma \ref{le:lemma14}. The proofs of Lemmas \ref{le:lemma3}-\ref{le:lemma5} are built on Lemmas \ref{le:lemma6}-\ref{le:lemma9}.
Lemmas~\ref{le:lemma6}-\ref{le:lemma7} give estimates for trajectory from start to $P_2$ ($c_2\lambda \,<\,s\,<\, \lambda$ in Region 2). Lemma~\ref{le:lemma8} gives estimate of the behaviour between $P_2$ and $P_4$ ($s\,<\,c_2\lambda$) and Lemma \ref{le:lemma9} for the behaviour between $P_4$ and $P_5$ ($c_2\lambda \,<\,s\,<\, c_1\lambda$ in Region 3).
Lemmas \ref{le:lemma6}-\ref{le:lemma9} use more new lemmas about which we inform later.
The section ends with the proof of Statement 3.

Before we start with the proofs of the Statements and Lemma \ref{le:lemma3} we introduce
Lemmas \ref{le:lemma4} and \ref{le:lemma5}.
Lemma \ref{le:lemma5} can be seen as corollary from Lemma \ref{le:lemma4}.
The proof of Lemma \ref{le:lemma3} is very similar to proof of Lemma \ref{le:lemma4}.
We wish to formulate the most general upper estimate for $x_6$ in Lemma~\ref{le:lemma4}.
We find such an estimate in the case $T$ intersects $s\,=\,c_2\lambda$ before escaping
Region~2 using auxiliary estimates for $x_1, \, x_2, \, x_4$ and $x_5$. For the estimate we need some notations and assumptions.

We introduce the following notations
\begin{equation}
 C_1\,=\,\left(1-c_1+ \ln \left(c_1\right)\right)\lambda,\qquad C_2\,=\,\left(c_1-c_2 +\ln \left(\frac{c_2}{c_1}\right)\right)\lambda,
\label{c12}
\end{equation}
and notice that $C_1,\, C_2\,<\,0$.
Moreover, we let
\begin{equation}
H_0\,=\,a+\lambda,\qquad H_1\,=\,a+c_1\lambda,\qquad H_2\,=\,a+c_2\lambda.
\label{h012}
\end{equation}
Next, we assume that
\begin{equation}
C_1 \,>\, - \left(\sqrt{x_0} -\sqrt{H_0}\right)^2,
\label{ass1}
\end{equation}
and defined the function $Q_1$ by
\begin{equation}\label{eq:Q1-eq}
Q_1\left(x\right) \,=\, x^2 - \left(H_0+C_1+x_0\right) x + H_0 x_0.
\end{equation}
If \eqref{ass1} is satisfied,
then $Q_1$ has a unique root $x^+_1$ in the interval $\left(\sqrt{H_0 x_0}, x_0\right)$.
Similarly, we assume that
\begin{equation}
C_2 \,>\, - \left(\sqrt{x^+_1} -\sqrt{H_1}\right)^2,
\label{ass2}
\end{equation}
and define the function $Q_2$ by
\begin{equation}\label{eq:Q2-eq}
Q_2\left(x\right)\,=\,x^2 -\left(H_1+C_2+x^+_1\right) x + H_1 x^+_1.
\end{equation}
Again, if \eqref{ass2} is satisfied we note that then $Q_2$ has a unique root $x^+_2$ in the interval $\left(\sqrt{H_1 x^ +_1}, x^+_1\right)$.
We also introduce a function $\theta$ and a number $\tilde{C}$ by
\begin{equation}
\theta \left(x\right)\,=\, x - H_2 \ln \left(x\right), \qquad \tilde{C}\,=\,\theta \left(x^+_2\right).
\label{thC}
\end{equation}
We make one more assumption
%
\begin{equation}
k \,=\, H_2\, e^\frac{\tilde C}{H_2} \,>\, 4.
\label{ass3}
\end{equation}
Also the following notations are needed
\begin{equation}
\hat z \,=\,\frac{1}{\sqrt{k^ 2-4k}},\qquad  M\,=\,\frac{c_2}{c_1} \qquad \text{and} \qquad A\,=\,\frac{\left(1-c_1\right)\lambda}{H_1}.
\label{zma}
\end{equation}
With these assumptions and notations we can formulate an upper estimate for $x_6$.

\begin{lemma}\label{le:lemma4}
Suppose assumptions \eqref{ass1},  \eqref{ass2} and \eqref{ass3} are satisfied.
Then the trajectory $T$ intersects $s\,=\,c_2\lambda$ before escaping Region 2 and for $x_6$ we have the estimate
\begin{equation*}
x_6 \,<\, \left(1+\hat z\right) M^A e^{-\frac{\tilde C}{H_2}}.
\end{equation*}
\end{lemma}

From the definition of $x^+_1$ and $x^+_2$,
we obtain,
for $i \,=\, 1,2$,
that $x^+_i \,>\, x^*_i$ where $x^*_i$ is the value of $x^+_i$ for $a \,=\, 0.1,\, \lambda \,=\, 0.1$ and $x_0 \,=\, 1$.
This will give us a new estimate as a corollary which we call Lemma~\ref{le:lemma5}.
To formulate the lemma we need the notation
\begin{equation*}
D^*_1\,=\,1-\frac{H^*_0}{x^*_1}, \qquad D^*_2\,=\,1-\frac{H^*_1}{x^*_2},
\end{equation*}
where $H^*_i$ are the values of $H_i$, $i \,=\, 0,1$, when $a \,=\, \lambda \,=\, 0.1$,
that is, $H_0^* \,=\, 0.2$ and
$H_1^* \,=\, \left(1 + c_1\right) \cdot 0.1$.
With these notations we can formulate next lemma.

\begin{lemma}\label{le:lemma5}
Suppose assumptions \eqref{ass1},  \eqref{ass2} and \eqref{ass3} are satisfied.
Then the trajectory $T$ intersects $s\,=\,c_2\lambda$ before escaping Region 2 and for $x_6$ we have the estimate
\begin{equation*}
x_6 \,<\, \left(1+\hat z\right) M^A e^{-\frac{\theta \left(\tilde x_2\right)}{H_2}},
\end{equation*}
where
\begin{equation*}
\tilde x_2\,=\,x_0 + \frac{C_1}{D^*_1} + \frac{C_2}{D^*_2}.
\end{equation*}
\end{lemma}

Because $C_i, \, i\,=\,1,2$ depend only on $\lambda$,
$D^*_i$ only on $\lambda$ and $c_1$,
$\tilde x_ 2$ does not depend on $a$, only on $\lambda,\, c_1$ and $x_0$.
If we choose $c_1\,=\,e^{-2},\,c_2\,=\,e^{-4}$ we are able to prove that
assumptions \eqref{ass1}, \eqref{ass2} and \eqref{ass3} are satisfied and get Lemma \ref{le:lemma3}.

Lemma \ref{le:lemma4} is based on Lemmas \ref{le:lemma6}, \ref{le:lemma8} and \ref{le:lemma9}.
We now give these lemmas and also Lemma \ref{le:lemma7} needed for Lemma \ref{le:lemma3}.
Lemma \ref{le:lemma7} can be seen as a corollary of Lemma \ref{le:lemma6}.

\begin{lemma}\label{le:lemma6}
Suppose assumptions \eqref{ass1} and \eqref{ass2} are satisfied.
Then the trajectory $T$ intersects $s\,=\,c_i\lambda, \, i\,=\,1,2$,
before escaping Region 2 at points $P_1\,=\,\left(x_1,c_1\lambda \right)$ and  $P_2\,=\,\left(x_2,c_2\lambda \right)$, where
%
\begin{equation}
x_1 \,>\, x_0 +\frac{C_1}{1-\frac{H_0}{x^+_1}} \qquad \text{and} \qquad
x_2 \,>\, x_0 +\frac{C_1}{1-\frac{H_0}{x^+_1}} +\frac{C_2}{1-\frac{H_1}{x^+_2}}.
\label{e2}
\end{equation}
Moreover, if $x_i^*$ are the values for $x_i^+$, $i\,=\,1,2$,
when $a\,=\,\lambda \,=\,0.1$ and $x_0\,=\,1$,
then the inequalities in \eqref{e2} remain valid for all $a,\lambda\,<\,0.1$
if $x_i^+$ are replaced by $x_i^*$.
\end{lemma}

Using Lemma \ref{le:lemma6} for special values of $c_i$ after calculating some quantities we get a corollary.

\begin{lemma}\label{le:lemma7}
Suppose that $c_1\,=\,e^{-2}$ and $c_2\,=\,e^{-4}$.
Then the trajectory $T$ intersects $s\,=\,c_2\lambda$ before escaping Region~2 at point $P_2\,=\,\left(x_2,c_2\lambda \right)$, where
\begin{equation*}
x_2 \,>\, x_0 - 3.8 \lambda.
\end{equation*}
\end{lemma}

\begin{lemma}\label{le:lemma8}
Let $T^*$ be a trajectory of system \eqref{pp} under conditions \eqref{pk} with initial conditions
$x\left(0\right)\,=\,u, \, s\left(0\right)\,=\,\lambda^* \,<\, \lambda, \, u\,>\,2.5H^*$ and $H^*\,=\,a+\lambda^*\,>\,h\left(\lambda^*\right)$.
Then the trajectory $T^*$ next time intersects $s\,=\,\lambda^*$ at a point $P\,=\,\left(v,\lambda^*\right)$ where
\begin{equation*}
v\,<\,\left(1+\hat z\right)e^{-\frac{\tilde\theta\left(u\right)}{H^*}},
\end{equation*}
and where
\begin{equation*}
\tilde\theta \left(u\right)\,=\,u-H^* \ln \left(u\right),\qquad
\hat z \,=\,\frac{1}{\sqrt{\tilde k^ 2-4\tilde k}} \qquad \text{and} \qquad
\tilde k\,=\,\frac{H^*}{u} e^\frac{u}{H^*}.
\end{equation*}
\end{lemma}

\begin{lemma}\label{le:lemma9}
Let $T^*$ be a trajectory of system \eqref{pp} under conditions \eqref{pk}
with initial conditions $x\left(0\right)\,=\,u, \, s\left(0\right)\,=\,c_2\lambda, \,
u\,<\,h\left(c_2\lambda\right)$. The trajectory $T^*$ next time intersects $s\,=\,c_1\lambda$ at a point $P_v\,=\,\left(v,c_1\lambda\right)$ where
\begin{equation*}
v\,<\,u M^A.
\end{equation*}
\end{lemma}

We now proceed to prove these lemmas.
We start with Lemma \ref{le:lemma6} and \ref{le:lemma7}.
The proof of Lemma \ref{le:lemma6} is based on Lemmas \ref{le:lemma10} and \ref{le:lemma11} which we give here,
before the proofs of Lemma \ref{le:lemma8} and \ref{le:lemma9}.
We consider a trajectory $T_2$ of system \eqref{pp} under conditions \eqref{pk}  with initial condition $x\left(0\right)\,=\,u\,>\,h\left(c_1^*\lambda\right), \, s\left(0\right)\,=\, c_1^*\lambda,\,
0\,<\,c_1^*\,\leq\, 1$.
Let $c_2^*$ be a number less than $c_1^*$.

\noindent
We introduce the quantities $C$ and $H$ and the function $R$ by
\begin{equation*}
C\,=\,\left(c_1^*-c_2^* + \ln \left(\frac{c_2^*}{c_1^*}\right)\right)\lambda, \qquad
H\,=\,a+c_1^*\lambda,\qquad R\left(x\right)\,=\,x-u- H\ln \left(\frac{x}{u}\right).
\end{equation*}
We are interested in whether $T_2$ intersects $s\,=\,c_2^*\lambda$ before escaping Region 2. We are also interested in a lower estimate for the $x$-value of such an intersection. Lemma \ref{le:lemma10} gives an answer to these questions and Lemma \ref{le:lemma11} gives a more explicit estimate.

\begin{lemma}\label{le:lemma10}
If the equation $R\left(x\right)\,=\,C$ has a solution $x\,=\,\bar x$, $H \,<\, \bar x \,<\, u$, then the trajectory $T_2$ intersects $s\,=\,c_2^*\lambda$ before escaping Region 2 at a point $\tilde P \,=\,\left(\tilde x, c_2^*\lambda\right)$, where $\tilde x \,>\, \bar x$.
\end{lemma}

Suppose $H, C$ and $u$ satisfy the following assumptions
\begin{equation}
 0\,<\,H\,<\,0.2, \qquad H\,<\,u, \qquad -\left(\sqrt{u}-\sqrt{H}\right)^2 \,<\, C \,<\,0,
\label{assC}
\end{equation}
and define
\begin{equation*}
\tilde R\left(x\right) \,=\,\left( 1-\frac{H}{x}\right) \left(x-u\right).
\end{equation*}
%
Then $Q\left(x\right) \,=\, \left(\tilde R\left(x\right)-C\right)x \,=\, x^2-\left(H+C+u\right)x + Hu \,=\, 0$ has a unique root
$x_+ \in \left(\sqrt{Hu}, u\right)$ and the following holds.

\begin{lemma}\label{le:lemma11}
If \eqref{assC} is satisfied, then equation $R\left(x\right)\,=\,C$ has exactly one solution $\bar x$ for $x\,>\,H$ and
%
\begin{equation}
\bar x \,>\, x_+ \,=\, u+ \frac{C}{1-\frac{H}{x_+}}.
\label{l2b1}
\end{equation}
Further if $H\,<\,H_m$ and $x_m$ is the root of $Q$ between $\sqrt{Hu}$ and $u$ for $H\,=\,H_m$ and $u\,=\,1$, then
\begin{equation}
\bar x \,>\, u + \frac{C}{1-\frac{H_m}{x_m}}.
\label{l2b2}
\end{equation}
\end{lemma}

\begin{figure}[h]
\begin{center}
\includegraphics[scale = 0.5]{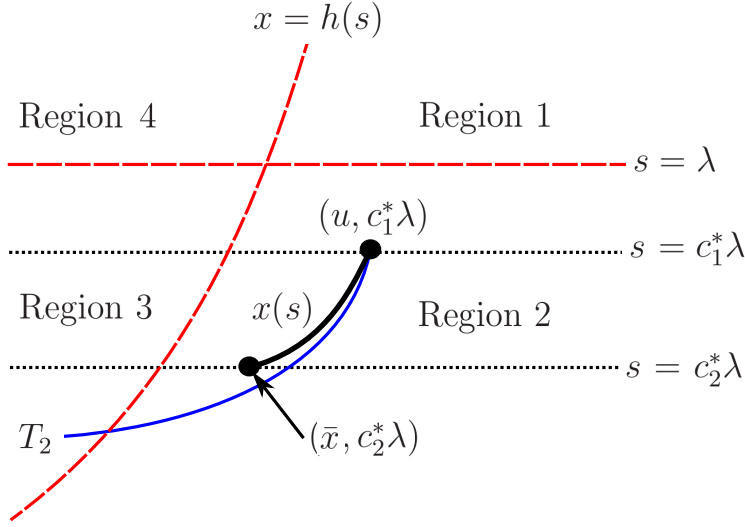}
\caption{Geometry in the proof of Lemma \ref{le:lemma10}.}
\label{fig:Lemma10}
\end{center}
\end{figure}

\noindent
{\it Proof of Lemma \ref{le:lemma10}}.
We notice that the equation $R\left(x\right)\,=\,C$ is equivalent to
$U\left(x,c_2^*\lambda\right) \,=\, U\left(u,c_1^*\lambda\right)$,
where
\begin{equation}\label{eq:U_def}
U\left(x,s\right) \,=\, x - H \ln\left(x\right) + s - \lambda \ln\left(s\right).
\end{equation}
Because the equation has a solution $\bar x$, $H \,<\, \bar x \,<\, u$,
and $U$ is increasing in $x$ and decreasing in $s$ in Region~2 as long as $H \,<\, x$,
the equation
$U\left(x,s\right) \,=\, U\left(u,c_1^*\lambda\right)$
has a unique solution $x\left(s\right)$ for any $s$ between $c_2^*\lambda$ and $c_1^*\lambda$ and $x\left(s\right)$ is
%
%
increasing in $s$ and $x\left(c_2^*\lambda\right) \,=\, \bar x$ and $x\left(c_1^*\lambda\right) \,=\, u$,
see Figure \ref{fig:Lemma10}.
Derivation with respect to time gives $U'\left(x,s\right) \,=\, \left(h\left(s\right)-H\right)\left(s-\lambda\right)\,>\,0$ in Region~2.
Thus, because $U\left(x,s\right)$ increases in $x$,
the trajectory $T_2$ will remain in the region defined by $x\,>\,x\left(s\right)$ until it intersects $s\,=\,c_2^*\lambda$ at $\tilde P$ with $\tilde x\,>\, \bar x$.
On trajectory part $U\left(x,s\right)\,>\,U\left(x\left(s\right),s\right)\,=\,U\left(u,c_1^*\lambda\right)$. $\hfill\Box$\\

\noindent
{\it Proof of Lemma \ref{le:lemma11}}.
We now use the auxiliary function $\tilde R\left(x\right)\,=\,\left( 1-\frac{H}{x}\right) \left(x-u\right)$.
From $\ln \left(\frac{x}{u}\right) \,>\, \frac{1}{x} \left(x-u\right)$ it follows that $R\left(x\right)\,<\,\tilde R\left(x\right)$ for $x \,<\, u$.
Equation $\tilde R\left(x\right) \,=\, C$ has a unique solution $x_+$ in $\left(\sqrt{Hu},u\right)$ when
\eqref{assC} is satisfied.
($\tilde R$ has a global minimum $ \tilde R\left(\sqrt{Hu}\right)\,=\, -\left(\sqrt{u}-\sqrt{H}\right)^2 $ and $ \tilde R \left(u\right)\,=\,0$).
Because $\tilde R\left(x\right) \,=\, C$ is equivalent with $Q\left(x\right) \,=\, 0$, $x_+$ is also the greatest root of $Q$.
Equation $R\left(x\right) \,=\, C$ has a unique solution $\bar x$ in $\left(H,u\right)$ such that $x_+\,<\,\bar x\,<\,u$,
because $R\left(u\right) \,=\, 0$ and
%
%
$R\left(x_+\right) \,<\, \tilde R\left(x_+\right) \,=\, C$
and $R$ is growing for $x \,>\, H$.
Now we notice that $Q\left(x\right)\,=\,0$ is equivalent to
$$
x\,=\,u+\frac{C}{1-\frac{H}{x}},
$$
from which we get \eqref{l2b1}.

To prove the second inequality we note that the function $Q$ is increasing in $H$ for $x\,<\,u$ and decreasing in $u$ for $x\,>\,H$,
from which we conclude that $x_+$ is decreasing in $H$ and increasing in $u$ and,
therefore, $\frac{H_m}{x_m}\,>\,\frac{H}{x_+}$ which implies \eqref{l2b2}.
Thus, both inequalities of the lemma are proved and the proof is complete. $\hfill\Box$\\

We can now prove Lemma \ref{le:lemma6}.\\

\noindent
{\it Proof of Lemma \ref{le:lemma6}}. From Lemma \ref{le:lemma10} and Lemma \ref{le:lemma11} with
$c_1^*\,=\,1, \, c_2^*\,=\,c_1, C\,=\,C_1,\, H\,=\,H_0, \, u\,=\,x_0$
it follows that the trajectory $T$ intersects $s\,=\,c_1\lambda$ before escaping Region 2,
and that for this intersection the first inequality in Lemma \ref{le:lemma6} holds.
Using Lemma \ref{le:lemma10} and Lemma \ref{le:lemma11} once again, this time with
$c_1^*\,=\,c_1, \, c_2^*\,=\,c_2,\, C\,=\,C_2,\, H\,=\,H_1, \, u\,=\,x_1$,
we conclude that $T$ also intersects $s\,=\,c_2\lambda$ before escaping Region 2,
and that for this intersection the second inequality in Lemma \ref{le:lemma6} holds.
Indeed, to see that we can apply Lemma \ref{le:lemma11} here we observe that
$
u \,=\, x_1 \,>\, x_1^+ 
\,>\, \sqrt{H_0 x_0} \,>\, a + \lambda \,>\, H_1.
$
Finally, we notice that $H_0$ and $H_1$ take their maximal values for $a\,=\,\lambda\,=\,0.1$.
Thus, the possibility to replace $x_i^+$ by $x_i^*$ follows from inequality \eqref{l2b2}.$\hfill\Box$\\

\noindent
{\it Proof of Lemma \ref{le:lemma7}}. 
We intend to use Lemma \ref{le:lemma6}.
Let $c_1 \,=\, e^{-2}$ and $c_2 \,=\, e^{-4}$.
Using \eqref{c12} and \eqref{h012} we find
$C_1 \,>\, -1.135 \lambda$,
$C_2 \,>\, -1.883 \lambda$,
$H_0 \,<\, 0.2$ and
$H_1 \,<\, 0.1136$.
Equation \eqref{eq:Q1-eq} with $x_0 \,=\, 1$ yields $x_1^* \,>\, 0.851$ and
\eqref{eq:Q2-eq} with $x_1^+ \,=\, 0.851$ yields $x_2^* \,>\, 0.620$.
Using these estimates we obtain
$D_1 \,=\, 1 - \frac{H_0}{x_1^*} \,>\, 0.764, \, D_2\,=\, 1 - \frac{H_1}{x_2^*} \,>\, 0.816$ and
$\frac{C_1}{D_1} + \frac{C_2}{D_2} \,>\, -3.8 \lambda$.
Now, we note that the above estimates imply assumptions \eqref{ass1} and \eqref{ass2},
and Lemma \ref{le:lemma7} now follows by an application of Lemma \ref{le:lemma6}. $\hfill\Box$\\

We have now proved Lemmas \ref{le:lemma6} and \ref{le:lemma7} and proceed to the proof of Lemma~\ref{le:lemma8}.
Proof of Lemma \ref{le:lemma8} is based on Lemmas \ref{le:lemma12} and \ref{le:lemma13}, we now introduce.
We consider a trajectory $T^*$ of system \eqref{pp} under conditions \eqref{pk} with initial conditions
$
x\left(0\right)\,=\,u, \, s\left(0\right) \,=\, \lambda^*\,<\,\lambda, \, u \,>\, a + \lambda^*.
$
Suppose $P\,=\,\left(v,\lambda^*\right)$ is the next intersection with $s\,=\,\lambda^*$.
Let further
$$
\theta \left(x\right)\,=\,x- H\ln \left(x\right).
$$
The following lemma which will also be used for Statement 3,
gives estimates for $v$.

\begin{lemma}\label{le:lemma12}
Let $\hat x$ be the solution to $\theta \left(x\right)\,=\,\theta \left(u\right), \, x\,<\,H\,=\,a+\lambda^*$ and $\check{x}$  the solution to $\theta \left(x\right)\,=\,\theta \left(u\right), \, x\,<\,H\,=\,a$.
Then for the next intersection of trajectory $T^*$ with $s \,=\, \lambda^*$ at $P\,=\,\left(v,\lambda^*\right)$,
it holds that  $\check{x}\,<\,v\,<\,\hat{x}$.
\end{lemma}

Next lemma gives estimate for the equation in previous lemma.

\begin{lemma}\label{le:lemma13}
Suppose that $C$ is a number such that $C\,>\,\theta\left(H\right)$ and suppose
$$
k\,=\,He^{C/H}\,>\,4, \qquad 0\,<\,H\,<\,0.2.
$$
Then the equation
\begin{equation}
\theta \left(x\right)\,=\,C,\qquad x\,<\,H,
\label{eq13}
\end{equation}
has a unique solution $\bar x$ such that
\begin{equation*}
e^{-C/H} \,<\, \bar x\,<\, \left(1+\hat z\right)\, e^{-C/H},
\end{equation*}
where $\hat z \,=\,\frac{1}{\sqrt{k^2-4k}}$.
\end{lemma}

\begin{figure}[h]
\begin{center}
\includegraphics[scale = 0.5]{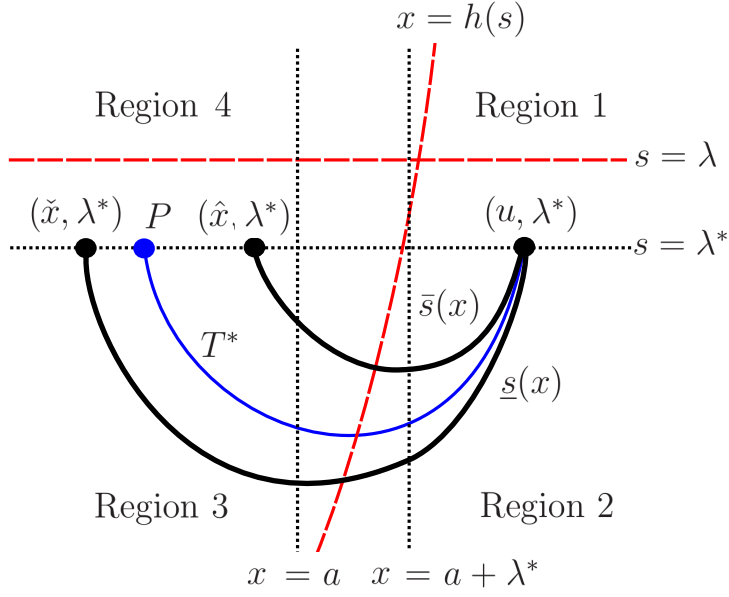}
\caption{Geometry in the proof of Lemma \ref{le:lemma12}.}
\label{fig:Lemma12}
\end{center}
\end{figure}

\noindent
{\it Proof of Lemma \ref{le:lemma12}}.
The trajectory $T^*$ escapes from Region 2 at a minimal $s$ to Region 3, where $s$ grows and after some time $T^*$ intersects $s\,=\,\lambda^*$ at $P\,=\,\left(v,\lambda^*\right)$, see Figure \ref{fig:Lemma12}.
As in the proof of Lemma \ref{le:lemma10} we will make use of the function $U$ defined in \eqref{eq:U_def}
to construct barriers for the trajectory $T^*$.
We note that $U$ is decreasing in $s$,
increasing in $x$ for $x \,>\, H$ and decreasing in $x$ for $x \,<\, H$.
Moreover, $\theta\left(x\right) \,=\, U\left(x,\lambda^*\right)$.

We first prove the upper bound $v \,<\, \hat{x}$.
Let $H\,=\,a+\lambda^*$ and let $\bar{s}\left(x\right)$ be the level curve to $U$ such that $U\left(x, \bar{s}\left(x\right)\right) \,=\, \theta\left(u\right)$.
The curve $\bar{s}\left(x\right)$ will have a minimum at $x \,=\, H \,=\, a+\lambda^*$ and
intersect $\lambda^*$ at $x \,=\, \hat x$ and also at $x \,=\, u$.
Observe that,
since $h\left(s\right) \,<\, h\left(\lambda^*\right) \,<\, \lambda^* + a \,=\, H$,
the derivative of $U$ with respect to time is positive: $U'\,=\,\left(h\left(s\right)-H\right)\left(s-\lambda\right) \,>\, 0$.
Therefore, the trajectory $T^*$ must stay below the curve $\bar{s}\left(x\right)$.
On trajectory $T^*$ we have $U\left(x,s\right) \,>\, U\left(x, \bar{s}\left(x\right)\right) \,=\, \theta\left(u\right) \,=\, \theta\left(\hat{x}\right)$.
Hence, recalling that $U$ is decreasing in $x$ for $x\,<\,H$, we have $v \,<\, \hat x$ and the upper bound follows.

The proof of the lower bound $\check{x} \,<\, v$ is similar.
Let $H \,=\, a$ and let $\underline{s}\left(x\right)$ be the level curve to $U$ such that $U\left(x, \underline{s}\left(x\right)\right) \,=\, \theta\left(u\right)$.
In this case, the derivative of $U$ with respect to time is negative,
and thus the trajectory $T^*$ must stay above the curve $\underline{s}\left(x\right)$.
On trajectory $T^*$ we have $U\left(x,s\right) \,<\, U\left(x, \underline{s}\left(x\right)\right) \,=\, \theta\left(u\right) \,=\, \theta\left(\check{x}\right)$,
and it follows also that $\check{x} \,<\, v$.  $\hfill\Box$\\

\noindent
{\it Proof of Lemma \ref{le:lemma13}}.
It is clear that equation \eqref{eq13} must have a solution because
$\theta \left(H\right) \,<\,C$ and $\theta \left(x\right) \to \infty$ for $x\to 0_+$.
The solution is unique because $\theta$ is decreasing for $x\,<\,H$.
Moreover, since $e^{-C/H} \,<\, H$ and $C \,<\, \theta \left(e^{-C/H}\right)$,
a solution $\bar x$ of $\theta \left(x\right) \,=\, C$ must satisfy $\bar x \,>\, e^{-C/H}$,
which proves the first inequality in Lemma \ref{le:lemma13}.

Substitution of $x\,=\,\left(1+z\right)e^{-C/H}$ into  $\theta \left(x\right)$  gives
$\bar\theta \left(z\right) +C$  where $\bar\theta \left(z\right) \,=\, \left(1+z\right)e^{-C/H} - H\ln \left(1+z\right)$.
Thus equation \eqref{eq13} is equivalent to $\bar\theta \left(z\right) \,=\,0$.
Let $\bar z$ be the $z$-value corresponding to the solution $\bar x$
($\bar x\,=\,\left(1+\bar z\right)e^{-C/H}$).
For $z\,=\,0$ we get $\theta \left(x\right)\,=\,x+C\,>\,C$, so clearly $\bar z \,>\,0$.
We wish to find an upper estimate for $\bar z$.
From $\ln \left(1+z\right)\,>\, \frac{z}{1+z}$ it follows that $\bar \theta \left(z\right)\,<\,
\left(1+z\right)e^{-C/H} - \frac{Hz}{1+z} \,=\, \frac{ e^{-C/H}}{1+z} \tilde \theta \left(z\right)$,
where $\tilde\theta \left(z\right)\,=\,\left(1+z\right)^ 2 -kz,\, k\,=\,H\, e^{C/H}$.
The function $\tilde \theta\left(z\right)$ has two roots because $k\,>\,4$.
We denote the smallest one by $\tilde z$.
Clearly  $\tilde z \,<\,\frac{k-2}{2}$,
and $1+\tilde z \,<\, k/2 \,<\,k$ which is equivalent to $\tilde x \,=\, \left(1+\tilde z\right)e^{-C/H} \,<\, H$.
Because $\bar \theta \left(\tilde z\right)\,<\,\tilde \theta \left(\tilde z\right)\,=\,0$ and $\bar\theta$ is decreasing in $\left(0,\tilde z\right)$ we must have $\bar z \,<\, \tilde z$.
 Using the assumption $k \,>\, 4$ and the mean value theorem,
we get an estimate for $\tilde z$:
$$
\tilde z \,=\,\frac{k-2}{2}\left( 1- \sqrt{1-\frac{4}{\left(k-2\right)^ 2}}\right)\,<\,
\frac{k-2}{2}\cdot \frac{4}{\left(k-2\right)^ 2}\cdot \frac{1}{2\sqrt{1-\frac{4}{\left(k-2\right)^ 2}}} \,=\, \frac{1}{\sqrt{k^ 2-4k}}\,=\, \hat z.
$$
Now we conclude $0\,<\,\bar z\,<\,\hat z$ and thereby $e^{-C/H} \,<\, \bar x \,<\,\left(1+\hat z\right)e^{-C/H}$ and the lemma is proved. $\hfill\Box$\\

We are now ready with the proofs of Lemma \ref{le:lemma12} and \ref{le:lemma13} and can use them for proving Lemma \ref{le:lemma8}.\\

\noindent
{\it Proof of Lemma \ref{le:lemma8}}.
The result follows from Lemma \ref{le:lemma12} and Lemma \ref{le:lemma13} by taking $H \,=\, H^*$ and $C \,=\, \tilde \theta\left(u\right)$.
We observe that we will have $\tilde k \,>\,4$ because $u\,>\,2.5H^*$. $\hfill\Box$\\

\noindent
Now only Lemma \ref{le:lemma9} is left to be proved in order to give the proofs of Lemmas~\ref{le:lemma3}-\ref{le:lemma5}.\\

\noindent
{\it Proof of Lemma \ref{le:lemma9}}.
The part of the trajectory between $P_u\,=\,\left(u,c_2\lambda\right)$ and $P_v$ is in Region 3,
where $s'\,>\,0\,>\,x'$ and moreover $s\,<\,c_1\lambda$.
There we get the following inequalities
\begin{equation*}
\frac{dx}{ds} \,<\, \frac{\max x'}{\max s'} \,=\,
\frac{\left(c_1-1\right)\lambda x}{\left(h\left(c_1\lambda\right)-x\right)s}
 \,<\,\frac{\left(c_1-1\right)\lambda x}{\left(a+c_1\lambda\right)s}.
\end{equation*}
Integrating and using $u \,<\, c_2 \lambda$ we get
$$
\int_u^ v \frac{1}{x} dx \,<\, \int_{c_2\lambda}^{c_1\lambda}
\frac{\left(c_1-1\right)\lambda}{a+c_1\lambda} \cdot \frac{1}{s} ds,
$$
and, by using the notation in \eqref{zma} we have
$
\ln\left(v\right) - \ln\left(u\right) \,<\, - A \ln\left(M^{-1}\right)
$
from which Lemma \ref{le:lemma9} follows. $\hfill\Box$ \\

We have now finished the proofs of all auxiliary results needed for Lemmas \ref{le:lemma3}-\ref{le:lemma5} and we will now continue by proving these lemmas.\\

\noindent
{\it Proof of Lemma \ref{le:lemma4}}.
From Lemma \ref{le:lemma6} follows that the trajectory $T$ intersects $s\,=\,c_2\lambda$
before escaping Region 2 at a point $P_2\,=\,\left(x_2,c_2\lambda\right)$ where $x_2 \,>\, x_2^+$.
From
Lemma \ref{le:lemma7} it follows that $x_2^+ \,>\, x_0 - 3.8 \lambda \,>\, 0.6$, and,
therefore, we can apply Lemma \ref{le:lemma8} with $H^*\,=\,H_2 \,=\, a + c_2 \lambda$.
In particular, from Lemma \ref{le:lemma8} with $H^*\,=\,H_2$ and $\lambda^* \,=\, c_2\lambda$
it follows that the trajectory with initial condition
$x\left(0\right)\,=\,x_2^+,\,s\left(0\right)\,=\,c_2\lambda\,=\,\lambda^*$ next time intersects $s\,=\,c_2\lambda$ at a point
$\tilde P_4\,=\,\left(\tilde x_4,c_2\lambda\right)$
where $\tilde x_4\,<\, \left(1+\hat z\right) e^{-\frac{\theta\left(x_2^ +\right)}{H_2}}$. Thus trajectory $T$ intersects
$s\,=\,c_2\lambda$ at a point $P_4\,=\,\left(x_4,c_2\lambda\right)$, where $x_4\,<\,\tilde x_4$.
From Lemma \ref{le:lemma9} follows that a trajectory with initial condition $x\left(0\right)\,=\,\tilde x_4,\, s\left(0\right)\,=\,c_2\lambda$ next time intersects
$s\,=\,c_1\lambda$ at a point $\tilde P_5\,=\,\left(\tilde x_5,c_1\lambda\right)$, where $\tilde x_5\,<\,\tilde x_4 M^A$.
Thus trajectory $T$ intersects  $s\,=\,c_1\lambda$ next time at a point $P_5\,=\,\left(x_5,c_1\lambda\right)$, where $x_5\,<\,\tilde x_5$.
Finally, because at $P_6\,=\,\left(x_6,\lambda\right)$
(next intersection of $T$ with $s\,=\,\lambda$),
$x_6\,<\,x_5$ we get
$$
x_6\,<\,x_5\,<\,\tilde x_5\,<\,\tilde x_4 M^A \,<\,\left(1+\hat z\right) M^A e^{-\frac{\theta\left(x_2^+\right)}{H_2}}
$$
The proof of Lemma \ref{le:lemma4} is complete. $\hfill\Box$\\

\noindent
{\it Proof of Lemma \ref{le:lemma5}}.
The proof is analogous to the proof of Lemma \ref{le:lemma4}.
We only use Lemma \ref{le:lemma6} so that we replace $x_2^+$ by $x_2^*$ and modify it by taking as $H_0$ and $H_1$ the values they get for $a\,=\,\lambda \,=\,0.1$. $\hfill\Box$\\

\noindent
{\it Proof of Lemma \ref{le:lemma3}}.
The proof is analogous to proof of Lemma \ref{le:lemma4}, we only use Lemma \ref{le:lemma7} instead of Lemma \ref{le:lemma6}.
In particular, from Lemma \ref{le:lemma7} it follows that the trajectory $T$ intersects $s \,=\, e^{-4}\lambda$
before escaping Region 2 at a point $P_2\,=\,\left(x_2, e^{-4}\lambda\right)$, where
$$
x_2 \,>\, \tilde x_2 \,=\, x_0 - 3.8 \lambda.
$$
We now use Lemma \ref{le:lemma8} with $H^* \,=\, a + \lambda^*$, $\lambda^* \,=\, e^{-4}\lambda$ and $u \,=\, \tilde{x_2}$ to obtain
$
x_4 \,<\, \tilde x_4 \,=\, \left(1 + \hat{z}\right) e^{-\frac{ \tilde{\theta} \left(\tilde{x}_2\right) }{a + e^{-4}\lambda}}.
$
To estimate $\hat{z}$, we carefully observe that the largest $\hat{z}$ is obtianed by setting
$a \,=\, \lambda \,=\, 0.1$ and $x_0 \,=\, 1$.
Indeed, we obtian $\hat{z} \,<\, 0.015$ and so
$$
x_4 \,<\, \tilde x_4 \,=\, 1.015 e^{-\frac{ \tilde{\theta} \left(\tilde{x}_2\right) }{a + e^{-4}\lambda}}.
$$
From Lemma \ref{le:lemma9} follows that a trajectory with initial condition $x\left(0\right)\,=\,\tilde x_4,\, s\left(0\right)\,=\,e^{-4} \lambda$ next time intersects $s\,=\,e^{-2}\lambda$ at a point $\tilde P_5\,=\,\left(\tilde x_5, e^{-2}\lambda\right)$,
where $\tilde x_5 \,<\, \tilde x_4 M^A \,=\, \tilde x_4 e^{-2 A}$.
Thus trajectory $T$ intersects  $s \,=\, e^{-2} \lambda$ next time at a point $P_5 \,=\, \left(x_5, e^{-2}\lambda\right)$,
where $x_5 \,<\, \tilde x_5$.
Finally, because at $P_6\,=\,\left(x_6,\lambda\right)$ we have $x_6 \,<\, x_5$, we get
$$
x_6\,<\,x_5\,<\,\tilde x_5\,<\,\tilde x_4 e^{-2 A} \,<\, 1.015 e^{-2 A - \frac{ \tilde{\theta} \left(\tilde{x}_2\right) }{a + e^{-4}\lambda}},
$$
which proves Lemma \ref{le:lemma3}. $\hfill\Box$\\

In order to prove Statement 2 we need one more lemma.
The proof of it follows by using Lemma \ref{le:lemma13} with
$C \,=\, \theta\left(u\right)$ and $H \,=\, a$,
but it can also be proved shortly directly.

\begin{lemma}\label{le:lemma14}
Equation $\theta \left(x\right)\,=\,\theta \left(u\right), \, u\,>\,1$, where
$\theta \left(x\right)\,=\,x- a\ln \left(x\right), \, a\,<\,0.1$, has a unique solution $\bar x$ in $\left(0,a\right)$ and $\bar x \,>\, e^{-\frac{u}{a}} \,=\, \check x$.
\end{lemma}

\noindent
{\it Proof}.
We first note that $\theta $ is decreasing in $\left(0,a\right)$
and that $\theta$ has its global minimum at $a$.
Moreover, $\theta\left(u\right) \,=\, u - a \ln\left(u\right) \,<\, u/a + e^{-u/a} \,=\, \theta\left(\check{x}\right)$.
Therefore, $ \theta \left(a\right) \,<\, \theta\left(u\right) \,<\, \theta \left(\check x\right)$ and thus there is a unique solution to
$\theta \left(x\right) \,=\, \theta \left(u\right)$ between $\check x$ and $a$. $\hfill\Box$\\

We have now finished the proofs of all auxiliary lemmas and will proceed to the proofs of our main results for this section;  Statement 2 and Statement 3.\\

\noindent
{\it Proof of Statement 2}.
Lemma \ref{le:lemma12} and Lemma \ref{le:lemma14} together give the lower estimate in Statement 2 if we use $\lambda^*\,=\,\lambda$ and $u\,=\,x_0$.

To prove the upper bound we first observe that from Statement 1 and Lemma \ref{le:lemma3} it follows that
$x_0 - 3.8 \lambda \,<\, \tilde x_2 \,<\, 1.6$, where $\tilde x_2$ is as defined in Lemma \ref{le:lemma3}.
Using this estimate we conclude, since $\frac{\ln \left(1.6\right)}{1.6} \,<\, 0.294$, that
$$
\theta\left(\tilde{x}_2\right) \,=\,
\tilde{x}_2 \left( 1 - H_2 \frac{\ln\left(\tilde{x}_2\right)}{\tilde{x}_2}\right) \,>\,
\tilde{x}_2\left(1 - 0.294 H_2\right) \,>\,
\left( x_0 - 3.8 \lambda\right)\left(1 - 0.294 H_2\right).
$$
From \eqref{ec3} in Lemma \ref{le:lemma3}, using that $ H_2 \,<\, H_1$, we get
\begin{align*}
\ln \left(x_6\right) &\,<\, -2A-\frac{\theta \left(\tilde x_2\right)}{H_2} + \ln \left(1.015\right) \\
		  &\,<\, \frac{-2\left(1-c_1\right)\lambda}{H_1} - \frac{\left(x_0-3.8\lambda\right)\left(1-0.294\, H_2\right) - H_2\ln\left(1.015\right)}{H_2} \\
		  &\,<\, \frac{-2\left(1-c_1\right)\lambda - \left(x_0-3.8\lambda\right)\left(1-0.294\, H_2\right) + \left(a+c_2\lambda\right) \ln \left(1.015\right)}{H_1},
\end{align*}
and hence, using that $1 \,<\, x_0$,
\begin{align*}
\ln \left(x_6\right) &\,<\, \frac{-x_0 + \left(3.8 - 2\left(1-c_1\right) + c_2\ln \left(1.015\right) + 0.294 c_2 x_0\right)\lambda  +  \left(0.294 x_0 + \ln \left(1.015\right)\right)a}{H_1} \\
		  &\,<\, \frac{-x_0 \left(1 - 2.1\lambda - 0.31 a\right)}{H_1}
		   \,<\, -\frac{x_0}{1.32 H_1}.
\end{align*}
%
The above inequality gives the upper estimate in Statement 2 and the proof is complete. $\hfill\Box$\\

\noindent
{\it Proof of Statement 3}.
We denote by $S_2$ the part of the trajectory $T$ between the initial point
$\left(x_0,\lambda\right)$ and $P_6 \,=\, \left(x_6, \lambda\right)$.
Let $U\left(x,s\right)\,=\,x - H \, \ln \left(x\right) + s- \lambda\, \ln \left(s\right)$.
We denote by $\check s$  the solution to
\begin{equation}
U\left(x_0,\lambda \right)\,=\,U\left(a,s\right),\qquad H\,=\,a.
\label{e3da}
\end{equation}
The function $U$, for $H\,=\,a$, is decreasing in $s$ ($s\,\leq\, \lambda$),
decreasing in $x$ for $x\,<\,a$ and increasing in $x$ for $x\,>\,a$.
Thus the solutions to $U\left(x_0,\lambda \right)\,=\,U\left(x,s\right)$,
$s \,\leq\, \lambda$,
form a curve $\check S_2$ given by $s\,=\,\check\sigma \left(x\right),\, x_6\,\leq\, x \,\leq\, x_0$,
where $\check\sigma$ has a minimum for $x\,=\,a$ and is increasing for $x\,>\,a$ and decreasing for $x\,<\,a$.
Differentiating $U$ with respect to time gives $U'\,=\,\left(h\left(\lambda \right)-H\right)\left(s-\lambda\right)\,<\,0$
for $s\,<\,\lambda$ and
$U\left(x,s\right) \,<\, U\left(x,\check \sigma \left(x\right)\right)$  and $s \,>\, \check\sigma \left(x\right)$ for $\left(x,s\right)$ on trajectory $T$.

We denote by $\hat s$ the solution to
%
\begin{equation}
U\left(x_0,\lambda \right)\,=\,U\left(a + \lambda,s\right),\qquad H\,=\,a+\lambda.
\label{e3db}
\end{equation}
Analogously we find that the solutions to $U\left(x_0,\lambda \right)\,=\,U\left(x,s\right)$,
$s \,\leq\, \lambda$ form a curve
$\hat S_2$ given by $s\,=\,\hat\sigma \left(x\right),\, x_6\,\leq\, x \,\leq\, x_0$,
where $\hat\sigma$ has a minimum for $x\,=\,a+\lambda$ and is increasing for $x\,>\,a+\lambda$ and decreasing for $x\,<\,a+\lambda$.
We now get $U'\,>\,0$ for $s\,<\,\lambda$ and hence
$U\left(x,s\right) \,>\, U\left(x,\hat \sigma \left(x\right)\right)$  and $s \,<\, \hat\sigma \left(x\right)$ for $\left(x,s\right)$ on trajectory $T$.
We conclude that $\hat S_2$ and $\check S_2$ together form a closed region and
$S_2$ is wholly inside this region. The $s$-values on $\hat S_2$ are greater than
the corresponding $s$-values for $T$ and the $s$-values on $\check S_2$ are less than the corresponding $s$-values for $T$, except at the coinciding endpoints of the curves $S_2$, $\check S_2$ and $\hat S_2$.

The minimum $s$-value on $\check S_2$ is $\check s$ and it must be less than the minimal $s$-value on $S_2$ and we get $\check s \,<\,s_3$. Analogously we get $\hat s \,>\, s_3$, where $\hat s$ is the minimum $s$-value on $\hat S_2$.
%
We will now find an estimate for the solution $\check{s}$ to equation \eqref{e3da}.
To do so we first note that equation \eqref{e3da} is equivalent to
\begin{equation*}
s-\lambda \ln \left(s\right) \,=\, x_0 - a\ln \left(x_0\right) + \rho \left(\lambda\right) -\rho\left(a\right) \,=\, L,
\end{equation*}
where $\rho\left(x\right)\,=\,x\left(1-\ln \left(x\right) \right)$.
Because $-\lambda \ln \left(s\right) \,<\, s-\lambda \ln \left(s\right)$ and, since $1 \,<\, x_0$,
$$
\hat L\,=\,x_0 +\rho \left(\lambda\right)-\rho \left(a\right) \,>\,L
$$
and since $s-\lambda \ln \left(s\right)$ decreases in $s$,
we get $\check s \,>\, e^{-\frac{\hat L}{\lambda}}$.
Moreover, from $\rho \left(a\right)\,>\,0$ and $x_0\,>\,1$ it follows that
$\hat L \,<\, x_0 +\rho \left(\lambda\right) \,<\, x_0 \left(1+\rho\left(\lambda\right)\right)$ and thus
$$
s_3 \,>\, \check s \,>\, e^{-\frac{\hat L}{\lambda}} \,>\, e^{-\frac{x_0}{\lambda \kappa_2}},
$$
where $\kappa_2 \,>\, \frac{1}{1 + \lambda \left(1 - \ln\left(\lambda\right)\right)}$.
This proves the lower estimate in Statement 3.

To prove the upper estimate in Statement 3 we will now find an estimate for the solution $\hat{s}$
to equation \eqref{e3db}.
To do so we first note that this equation is equivalent to
\begin{equation}
s-\lambda \ln \left(s\right) \,=\,x_0 - H \ln \left(x_0\right) + \rho \left(\lambda\right) -\rho\left(a+\lambda\right) \,=\, \check L.
\label{e3dc}
\end{equation}
To estimate the solution of \eqref{e3dc} we will make use of Lemma \ref{le:lemma13}.
In particular, Lemma \ref{le:lemma13} with $H \,=\, \lambda$ and $C \,=\, \check L$ gives
%
\begin{align}\label{eq:hat_s_bound_above}
\hat s \,<\, \left(1 + \hat z\right) e^{-\frac{\check{L}}{\lambda}}\qquad \text{where} \qquad \hat{z} \,=\, \frac{1}{\sqrt{k^2 - 4k}}.
\end{align}
Next, we find a lower estimate of $k$.
Using the inequality $\check L \,>\, 1 + \rho \left(\lambda\right) - \rho\left(a + \lambda\right) \,=\, \tilde L$,
we see that
$k \,=\, \lambda e^{\frac{\check L}{\lambda}} \,>\, \lambda e^{\frac{\tilde L}{\lambda}} \,=\, \tilde k$.
The derivative with respect to $a$ of $\tilde{k}$ is of the same sign as the derivative of
$-\rho \left(a+\lambda\right)$ with respect to $a$ which is negative.
$\tilde k$ can also be written in form
$$
\tilde k\,=\,e^{\frac{1}{\lambda}-\frac{a}{\lambda} +\left( 1+ \frac{a}{\lambda} \right) \ln \left(a+\lambda \right)}
$$
and the derivative of $\tilde{k}$ with respect to $\lambda$ can be seen to be
$$
\frac{\partial\tilde k}{\partial\lambda} \,=\, \frac{\tilde k}{\lambda^2} \left( \lambda + a-1 - a \ln \left(a+\lambda\right)\right)
$$
which again is negative for our values of $a$ and $\lambda$.
Thus, $\tilde k$ is greater than the value ($\,>\,324$) it takes of $a\,=\,\lambda\,=\,0.1$
and it follows that $\hat z\,<\, \frac{1}{\sqrt{\tilde k^2 -4\tilde k}} \,<\, 0.004$.
Using this estimate and \eqref{eq:hat_s_bound_above} we conclude that
$\hat s \,<\, 1.004 e^{-\frac{\check L}{\lambda}}$.

For trajectory $T$ we have $x_0 \,<\, 1.6$ and therefore $x_0 - H\ln \left(x_0\right) \,>\, x_0\left(1-0.294 \,H\right)$.
Moreover, $\rho \left(a+\lambda \right)-\rho \left(\lambda\right)\,<\,\rho \left(a\right)$ and we obtain
\begin{align*}
\ln \left(\hat s\right) \,<\, -\frac{x_0\left(1-0.294\, H\right) -\rho \left(a\right)}{\lambda} + \ln \left(1.004\right),
\end{align*}
and so
\begin{align*}
\ln \left(\hat s\right)
&\,<\, -\frac{x_0}{\lambda }\left( 1-0.294\, H -\rho \left(a\right) -\lambda \ln \left(1.004\right)\right) \\
&\,<\, -\frac{x_0}{\lambda }\left( 1 - 0.3\, \lambda -a\left(1.3 - \ln\left(a\right) \right) \right).
\end{align*}
%
Thus $s_3\,<\, \hat s \,<\, e^{-\frac{x_0}{\lambda\kappa_3}}$,
where $\kappa_3$ is as in Statement 3,
and the proof is complete. $\hfill\Box$


\setcounter{equation}{0} \setcounter{theorem}{0}

\section{Estimates in Region 4}
\label{sec:Region4}

We again consider a trajectory $T$ of system \eqref{pp} under conditions \eqref{pk}
with initial condition $x\left(0\right)\,=\,x_0\,>\,1, \, s\left(0\right)\,=\,\lambda$. We are interested in the behaviour of the trajectory in Region 4. The trajectory enters Region 4 at point
$P_6\,=\,\left(x_6,\lambda\right)$. We are interested in the next intersection of the trajectory with $s\,=\,s_7\,>\,0.5$ at point $P_7\,=\,\left(x_7,s_7\right)$ (if it occurs before escaping Region 4) and of the next  intersection with the isocline $x\,=\,h\left(s\right)$ at $P_8\,=\,\left(x_8,s_8\right)$, where $x_8\,=\,h\left(s_8\right)$. Lemma \ref{le:lemma3} from previous section gives an estimate for $x_6$ and we are able to show that for such $x_6$ the trajectory will intersect $s\,=\,0.8$ before escaping Region 4 and the escaping occurs at $P_8$, where $s_8\,>\,0.9$.

The main result is Statement 4 which is based on the following two lemmas.

\begin{lemma}\label{le:lemma15}
The trajectory $T$ after intersecting $s\,=\,\lambda$ next time always intersects $s\,=\,0.8$ at a point $P_7\,=\,\left(x_7,0.8\right)$,
where $x_7 \,<\, 0.012$, before escaping Region 4.
\end{lemma}

\begin{lemma}\label{le:lemma16}
If  trajectory $T$ after intersecting $s\,=\,\lambda$ next time intersects $s\,=\,0.8$ at a point $P_7\,=\,\left(x_7,0.8\right)$
where $x_7 \,<\, 0.012$ then it intersects the isocline $s'\,=\, 0$ next time for an $s$-value greater than $0.9.$
\end{lemma}

From these lemmas follows

\begin{statement}
Trajectory $T$ after intersecting $s\,=\,\lambda$ at $P_4$ escapes from Region 4 at  an $s$-value greater than 0.9.
\end{statement}

The trajectory in Region 4 is well estimated by $x\,=\,x_6 B$ for $s\,<\,0.8$, where
$B$ is defined in \eqref{be1234}.
For $s\,>\,0.8$ expression \eqref{5c}
gives a one-sided estimate for the trajectory while remaining in Region 4 and \eqref{5d}
gives an estimate for $s_8$ substituting $m\,=\,0.8$.

The proof of Lemma \ref{le:lemma16} is at the end of the section.
The proof of  Lemma~\ref{le:lemma15} is  based on some lemmas we provide here.
Lemma \ref{le:lemma15} uses Lemma \ref{le:lemma18} and Lemma \ref{le:lemma19}.
Lemma \ref{le:lemma18} gives us necessary conditions in form of inequalities for the trajectory to intersect $s\,=\,s_7$ before escaping Region 4 and an estimate for $x$-value at intersection point $P_7$.
Lemma~\ref{le:lemma19} tells us that we have to check the inequalities only for $a,\lambda\,=\,0.1$ to be sure they hold for all other parameters. Lemma \ref{le:lemma18} is based on Lemma \ref{le:lemma17} and Lemma~\ref{le:lemma3},
where Lemma \ref{le:lemma17} gives estimates in Region 4 and Lemma \ref{le:lemma3} takes care of estimates for trajectory in Regions 2 and 3. Lemma \ref{le:lemma17} is based on Lemma \ref{le:lemma20} and \ref{le:lemma21}.
Lemma \ref{le:lemma20} gives us estimates for trajectory in a part of Region 4 and Lemma \ref{le:lemma21} tells us that we need to check these estimates only for $s\,=\,\lambda$ and $s\,=\,s_7$ in order to be sure the trajectory will stay in the region.

We now give Lemmas \ref{le:lemma17}-\ref{le:lemma19}.
Lemma \ref{le:lemma19} can be proved directly,
but the proof of Lemma \ref{le:lemma17}, which is needed for proving Lemma \ref{le:lemma18},
needs more lemmas and is given later.

\begin{lemma}\label{le:lemma17}
Suppose $0\,<\,k\,<\,1$ and $s_7\,\geq\, 0.5$. If
\begin{equation}
x_6\,  K_7 \left(\frac{1}{a+\lambda}\right)^{\frac{1}{k}}\,<\,\left(1-k\right)\, h\left(s_7\right)
\label{e4c}
\end{equation}
where
$K_7 \,=\, \left( e^{\frac{\lambda}{s_7}} \frac{s_7+a}{1-s_7} \right)^{\frac{1}{k}}$
and
$x_6 \,<\, \left(1-k\right)\, h\left(\lambda\right)$,
then the trajectory $T$ intersects $s \,=\, s_7$ before escaping Region 4 and at the intersection $P_7\,=\,\left(x_7,s_7\right)$ the estimate
\begin{equation*}
x_7 \,<\,  x_6\, K_7 \left(\frac{1}{a+\lambda}\right)^{\frac{1}{k}}
\end{equation*}
is satisfied.
\end{lemma}

Lemma \ref{le:lemma3} gives estimate for $x_6$ and thus
from Lemma \ref{le:lemma3} and Lemma \ref{le:lemma17} we get a new statement.
For this we introduce the function $\eta\left(a, \lambda\right)$ to find out whether inequality \eqref{e4c} holds.
\begin{equation*}
\eta \left(a,\lambda\right) \,=\, K_7  \left(\frac{1}{a+\lambda}\right)^{\frac{1}{k}}1.015\, e^{-2A-\frac{\theta \left(\tilde x_2\right)}{H_2}}
\end{equation*}
(Notations from Lemma \ref{le:lemma3} are used here).

\begin{lemma}\label{le:lemma18}
Suppose $0\,<\,k\,<\,1$ and $s_7\,\geq\, 0.5$. If
\begin{equation}
\eta \left(a,\lambda \right)\,<\,\left(1-k\right)\, h\left(s_7\right)
\label{e4d}
\end{equation}
and $x_6\,<\,\left(1-k\right)\, h\left(\lambda\right)$, then the trajectory $T$ intersects $s \,=\, s_7$ before escaping Region 4 and at the intersection $P_7\,=\,\left(x_7,s_7\right)$ the estimate $x_7\,<\,\eta \left(a,\lambda \right)$ holds.
\end{lemma}


\noindent
{\it Proof}. 
The statement follows directly from Lemma \ref{le:lemma17} and Lemma \ref{le:lemma3}. $\hfill\Box$\\

The following lemma tells us that to prove that \eqref{e4d} holds for all $a, \lambda \in [0,0.1)$,
it is enough to prove the inequality for $a \,=\, \lambda \,=\, 0.1$ fixed in the left hand side,
i.e., $\eta\left(0.1, 0.1\right) \,<\, \left(1-k\right)\, h\left(s_7\right)$.

\begin{lemma}\label{le:lemma19}
The derivatives of $\eta\left(a,\lambda\right)$ with respect to $a$ and $\lambda$ are positive if $k\,\geq\, 0.9$.
\end{lemma}

\noindent
{\it Proof}.
Calculations give
$\frac{\partial\eta}{\partial\lambda} \,=\, \xi_\lambda \eta\left(a,\lambda\right)$,
where
\begin{equation*}
\xi_\lambda \,=\, \frac{1}{s_7} - \frac{1}{k\left(a+\lambda\right)}-2\frac{\partial A}{\partial \lambda} + 3.8 \left( \frac{1}{H_2} -\frac{1}{\tilde x_2}\right) + \frac{\tilde x_2 c_2}{H^2_2}
\end{equation*}
and
\begin{equation*}
\frac{\partial A}{\partial \lambda} \,=\, \frac{\left(1-c_1\right)a}{H_1^2}\,<\,\frac{0.9a}{H_1^2} \,<\, \frac{0.9}{H_1}.
\end{equation*}
For $\xi_\lambda$ we get the estimate
$$
\xi_\lambda \,>\,  -\frac{1}{k\left(a+\lambda\right)} - \frac{1.8}{H_1} + \frac{3.8}{H_2} - \frac{3.8}{0.6} \,>\, -\frac{3}{H_2} + \frac{3.8}{H_2}  - \frac{3.8}{0.6}\,>\, \frac{0.8}{H_2}  - \frac{3.8}{0.6} \,>\, 0,
$$
because $H_2\,<\,0.11$.
Thus, since $\eta\left(a,\lambda\right) \,>\, 0$ we conclude that $\eta\left(a,\lambda\right)$ is increasing in $\lambda$.
Calculations also give $\frac{\partial\eta}{\partial a} \,=\, \xi_a \eta\left(a,\lambda\right)$, where
\begin{equation*}
\xi_a \,=\, \frac{1}{k \left(s_7 + a\right)} - \frac{1}{k\left(a+\lambda\right)}+ 2\frac{\left(1-c_1\right)\lambda}{H_1^2} + \frac{\tilde x_2 }{H^2_2}.
\end{equation*}
For $\xi_a$ we get the estimate
$$
\xi_a \,>\,  -\frac{1}{k\left(a+\lambda\right)} + \frac{\tilde x_2 }{H^2_2} \,>\, -\frac{1}{k\, H_2}+ \frac{0.6 }{H^2_2} \,=\, \frac{1}{H_2}\left( \frac{0.6}{H_2} - \frac{1}{k}\right) \,>\, 0.
$$
Thus, $\eta\left(a,\lambda\right)$ is increasing also in $a$ and the proof of the Lemma \ref{le:lemma19} is complete. $\hfill\Box$ \\


\vskip0.3cm

We now proceed to the proof of
Lemma \ref{le:lemma17}, which follows from Lemma \ref{le:lemma20} and \ref{le:lemma21}.
These lemmas  we formulate now.
For the statements we need the following notations:
\begin{equation}
B\,=\,E_1\,  E_2\,  E_3\,  E_4,
\label{be1234}
\end{equation}
where
%
\begin{equation}
E_1\,=\,\left( \frac{s+a}{s} \right)^{k_2}, \qquad
E_2\,=\,\left( \frac{\lambda}{a+\lambda} \right)^{k_2}, \qquad
\label{e1234a}
\end{equation}
\begin{equation}
E_3\,=\,\left( \frac{\left(s+a\right)\left(1-\lambda\right)}{1-s} \right)^{k_3}, \qquad
E_4\,=\,\left( \frac{1}{a+\lambda} \right)^{k_3},
\label{e1234b}
\end{equation}
\begin{equation*}
k_2\,=\,\frac{\lambda}{a} \qquad \textrm{and} \qquad k_3\,=\,\frac{1-\lambda}{1+a}.
\end{equation*}

\begin{lemma}\label{le:lemma20}
Suppose that $x_6\,<\,\left(1-k\right)\, h\left(\lambda\right)$.
Then as long as the trajectory stays in the region determined by
$x\,<\,\left(1-k\right)\, h\left(s\right)$ the following estimates are valid:
\begin{equation}
x\,<\,x_6B^{\frac{1}{k}}
\label{e7a}
\end{equation}
and,
when $s \,\geq\, 0.5$,
\begin{equation}
x\,<\,x_6\, K \left( \frac{1}{a+\lambda} \right)^{\frac{1}{k}}\qquad \text{where}\qquad
K\,=\,\left( e^{\frac{\lambda}{s}} \cdot \frac{s+a}{1-s}\right)^{\frac{1}{k}}.
\label{e7b}
\end{equation}
\end{lemma}

\begin{lemma}\label{le:lemma21}
Suppose $x_6\,<\,\left(1-k\right)\, h\left(\lambda\right)$ and $x_6 \left( B_7 \right)^{\frac{1}{k}}\,<\, \left(1-k\right)\, h\left(s_7\right)$,
where $B_7$ is the value $B$ takes for $s\,=\,s_7\,\geq\, 0.5$. Then the trajectory $T$ intersects $x\,=\,\left(1-k\right)\, h\left(s\right)$ next time after $P_6$ for $s\,>\,s_7$ and is inside the region determined by $x\,<\,\left(1-k\right)\, h\left(s\right)$ before it intersects $s\,=\,s_7$.
\end{lemma}

\noindent
{\it Proof of Lemma \ref{le:lemma20}}.
When $x\,<\,\left(1-k\right)\, h\left(s\right)$ we have $s'\,>\,k h\left(s\right) s\,>\,0$ and $x'\,>\,0$. Thus we get the inequality
\begin{equation*}
\frac{dx}{ds} \,<\, \frac{\left(s-\lambda\right)s}{k h\left(s\right) s}.
\end{equation*}
Integrating gives
%
\begin{equation}\label{eq:F_def}
x\,<\, x_6 \left( \frac{F\left(s\right)}{F\left(\lambda\right)}\right)^{\frac{1}{k}}\qquad \text{where}\qquad
F\left(y\right)\,=\,\frac{\left(y+a\right)^{k_1}}{y^{k_2} \left(1-y\right)^{k_3}}
\end{equation}
and where $k_1\,=\,\frac{a+\lambda}{a\left(a+1\right)}$.
But
\begin{equation*}
\frac{F\left(s\right)}{F\left(\lambda\right)}\,=\,
\left( \frac{s+a}{\lambda +a}\right)^{k_1}
\left( \frac{\lambda}{s}\right)^{k_2}
\left( \frac{1-\lambda}{1-s}\right)^{k_3}
\end{equation*}
and because $k_1\,=\,k_2+k_3$ we get
%
\begin{equation*}
\frac{F\left(s\right)}{F\left(\lambda\right)}\,=\,
\left( \frac{s+a}{s}\right)^{k_2}
\left( \frac{\lambda}{\lambda +a}\right)^{k_2}
\left( \frac{\left(1-\lambda\right)\left(s+a\right)}{1-s}\right)^{k_3}
\left( \frac{1}{\lambda +a}\right)^{k_3}.
\end{equation*}
Choosing $B \,=\, \frac{F\left(s\right)}{F\left(\lambda\right)} \,=\, E_1 E_2 E_3 E_4$ we get estimate \eqref{e7a}.

To prove \eqref{e7b} we note that for $E_i, \, i\,=\,1,2,3,4$ we get the following estimates
%
$$
E_1\,=\, \left(\left( 1+ \frac{1}{\frac{s}{a}}\right) ^\frac{s}{a} \right)^\frac{\lambda}{s} \,<\, e^\frac{\lambda}{s},\qquad
E_2\,=\,\frac{1}{\left( 1+\frac{1}{k_2}\right)^{k_2}}\,<\,1
$$
$$
k_3\,<\,1,\qquad
E_3\,<\, \left(\frac{s+a}{1-s}\right)^{k_3} \,<\, \frac{s+a}{1-s},\qquad
E_4\,<\, \frac{1}{\lambda +a}.
$$
%
All these estimates together give \eqref{e7b}. $\hfill\Box$\\

\noindent
{\it Proof of Lemma \ref{le:lemma21}}.
We first claim that the assumptions in the lemma implies
\begin{equation}\label{eq:lemma21_claim}
\frac{B\left(s\right)^{\frac{1}{k}}}{h\left(s\right)} \,<\, \frac{1-k}{x_6} \qquad \text{for all} \qquad s \in [\lambda, s_7].
\end{equation}
Next, assume, by way of contradiction,
that the trajectory $T$ intersects the curve $x \,=\, \left(1-k\right)h\left(s\right)$ for some $s \in [\lambda, s_7]$.
Using claim \eqref{eq:lemma21_claim} we then obtain $x_6 B\left(s\right)^{\frac{1}{k}} \,<\, x$ for the point of intersection.
But from \eqref{e7a} in Lemma \ref{le:lemma20} it follows that $x \,<\, x_6 B\left(s\right)^{\frac{1}{k}}$ as long as $T$ stays in the region defined by
$x \,<\, \left(1-k\right)h\left(s\right)$. Using continuity this leads to a contradiction.
Hence, we conclude that the trajectory $T$ intersects $x\,=\,\left(1-k\right)\, h\left(s\right)$ next time, after $P_6$, for $s\,>\,s_7$
and $T$ is inside the region determined by $x\,<\,\left(1-k\right)\, h\left(s\right)$ before it intersects $s\,=\,s_7$.

To finish the proof of Lemma \ref{le:lemma21} it remains to prove that claim \eqref{eq:lemma21_claim} holds true.
To do so we observe that differentiating $G\left(s\right)\,=\,\frac{F\left(s\right)^\frac{1}{k}}{h\left(s\right)}$ with respect to $s$,
where $F\left(s\right)$ is given by \eqref{eq:F_def},
gives
$$
G'\left(s\right)\,=\,\frac{F\left(s\right)^\frac{1}{k}}{k\left(1-s\right)^2 s \left(s+a\right)^2} G^*\left(s\right),
$$
where
$$
G^*\left(s\right)\,=\,2ks^2 +\left(ak-k+1\right)s-\lambda.
$$
We conclude that $G\left(s\right)$ has a unique minimum between $s\,=\,\lambda$ and $s\,=\,0.5$,
when $a,\lambda \,<\, 0.1$ and $0\,<\,k\,<\,1$,
because
$$
\frac{9\lambda \left(a+2\lambda -1\right)}{10} \,=\, G^*\left(\lambda\right) \,<\, 0 \,<\, G^*\left(0.5\right) \,=\, \frac{9a-20\lambda +10}{20},
$$
and $G^*$ is increasing in $s$.
Thus,
the maximal value of $G$ in $\lbrack \lambda ,s_7\rbrack$ is either $G\left(\lambda\right)$ or $G\left(s_7\right)$.
Claim \eqref{eq:lemma21_claim} now follows since
$
\frac{G\left(s\right)}{F\left(\lambda\right)^{\frac{1}{k}}} \,=\, \frac{B\left(s\right)^\frac{1}{k}}{h\left(s\right)}
$
%
and the assumptions in the lemma equals
$$
\frac{B\left(\lambda\right)^\frac{1}{k}}{h\left(\lambda\right)} \,<\, \frac{1-k}{x_6} \qquad \textrm{and} \qquad \frac{B\left(s_7\right)^\frac{1}{k}}{h\left(s_7\right)} \,<\, \frac{1-k}{x_6}.
$$
The proof of Lemma \ref{le:lemma21} is complete. $\hfill\Box$\\


We are now ready with proofs of Lemmas \ref{le:lemma20} and \ref{le:lemma21} and can use them for getting
proofs of Lemma \ref{le:lemma17} and \ref{le:lemma18}.\\

\noindent
{\it Proof of Lemma \ref{le:lemma17}}. The proof follows from Lemma \ref{le:lemma20} and \ref{le:lemma21}.
Lemma~\ref{le:lemma21} tells that the trajectory will be inside the region $x\,<\,\left(1-k\right)\, h\left(s\right)$ and then Lemma \ref{le:lemma20} gives us the necessary estimates. $\hfill\Box$\\

Finally,
we are ready with all proofs of auxiliary results and can prove the main Lemmas \ref{le:lemma15} and \ref{le:lemma16} from which Statement 4 follows.\\

\noindent
{\it Proof of Lemma \ref{le:lemma15}}.
We choose $k\,=\,0.9$ and $s_7\,=\,0.8$ and calculate
$\eta \left(0.1,0.1\right) \,<\, 0.012 \,<\, \left(1-k\right)\, h\left(s_7\right)$
and then from Lemma \ref{le:lemma19} it follows that inequality \eqref{e4d} holds for all $a,\lambda \in [0,0.1)$.
Since $\eta \left(a,\lambda\right) \,<\, 0.012$ it follows that
$x_6\, K_7  \left(\frac{1}{a+\lambda}\right)^{\frac{1}{k}} \,<\, 0.012$
and because $K_7 \,>\, 1$ we also get,
using $k \,=\, 0.9$,
that $x_6 \,<\, 0.012 \cdot \left(a+\lambda\right) \,<\, \left(1-k\right)\, h\left(\lambda\right)$.
Lemma \ref{le:lemma15} now follows by an application of Lemma \ref{le:lemma18}. $\hfill\Box$\\

\noindent
{\it Proof of Lemma \ref{le:lemma16}}.
We consider trajectories of system \eqref{pp} in region
$$
E_{m,s_7} \,=\,\lbrace \left(x,s\right) \vert \, 0\,<\,x\,<\,m\left(1-s\right), s\,>\,s_7\rbrace ,
$$
where $m\,\leq\, s_7$.
Observe that
$$
m\left(1-s\right)\,\leq\, s_7\left(1-s\right)\,<\,s\left(1-s\right)\,<\,\left(s+a\right)\left(1-s\right)\,=\,h\left(s\right).
$$
In $E_{m,s_7}$ we get the estimates
%
\begin{equation}
s'\,>\,\left(m\left(1-s\right)-x\right)s, \quad x'\,<\,sx \qquad \text{and hence} \qquad \frac{ds}{dx} \,>\, \frac{m\left(1-s\right) - x}{x}.
\label{5a}
\end{equation}
Let us consider a trajectory with initial condition $x\left(0\right)\,=\,x_7,\, s\left(0\right)\,=\,s_7$, where $x_7\,<\,m\left(1-s_7\right)$.
Using \eqref{5a},
we conclude that as long as this trajectory remains in $E_{m,s_7}$,
it will be in the subregion bounded by the trajectory of the linear system
\begin{equation}
s'\,=\, m\left(1-s\right)-x,\quad x'\,=\,x,
\label{5b}
\end{equation}
with initial condition  $x\left(0\right)\,=\,x_7, s\left(0\right)\,=\,s_7$ and the lines $x\,=\,m\, \left(1-s\right)$ and $x\,=\,0$.
Solving system \eqref{5b} we find that the trajectory follows the curve
%
\begin{equation}
s\,=\,d\left( \frac{x_7}{x} \right)^m +1 -\frac{x}{1+m} \qquad \text{with} \qquad
d\,=\,s_7+\frac{x_7}{1+m} -1.
\label{5c}
\end{equation}
The trajectory leaves  $E_{m,s_7}$ when $x\,=\,m\, \left(1-s\right)$.
(Observe that then $s'\,=\,0$ for \eqref{5b}).
Substituting $x\,=\,m\, \left(1-s\right)$ into \eqref{5c} we get
$$
\frac{d x_7^m}{m^m\left(1-s\right)^m} +1-s -\frac{m}{m+1} \left(1-s\right) \,=\, 0,
$$
which is equivalent to
\begin{equation}
1-s \,=\, \frac{\left(-d\right)^\frac{1}{m+1} x_7^\frac{m}{m+1}}{m^\frac{m}{m+1}} \left(1+m\right)^\frac{1}{m+1}.
\label{5d}
\end{equation}
The above expression for $1-s$ increases with $x_7$,
for all $m \,\geq\, 0$, because
%
$$
\frac{\partial}{\partial x_7} \left(-d x_7^m\right) \,=\, - \left( \frac{x_7}{1+m} + m\left( s_7 +\frac{x_7}{m+1} -1\right) \right) x_7^{m-1} \,=\, x_7^{m-1} \left(m\left(1-s_7\right)-x_7\right) \,>\, 0.
$$
Thus, a lower boundary for the maximal $s$ can be calculated from \eqref{5d} for given $\left(x_7,s_7\right)$ choosing $m \,=\, s_7$.
Calculations show that if $s_7 \,=\, 0.8$ and $x_7 \,<\, 0.012$, then the maximal $s$ is greater than 0.9379.
We observe that system \eqref{5b} is not depending on $a$ or $\lambda$. Hence, the results are independent of these parameters.
The proof of Lemma \ref{le:lemma16} is complete.
$\hfill\Box$


\newpage

\setcounter{equation}{0} \setcounter{theorem}{0}

\section{Numerical results}
\label{sec:numerics}

Before comparing our analytical estimates to numerical simulations,
let us mention that to achieve accurate numerics of system \eqref{pp}
under assumption \eqref{pk} we recommend to transform the equations (e.g. log transformations)
to avoid variables taking on very small values.
Imposing linear approximations near the unstable equilibria at
$(x,s) = (0,0)$ and $(x,s) = (0,1)$ are also helpful.
Indeed, implementing MATLABs ode-solver ODE45 directly on system \eqref{pp}
may result in trajectories not satisfying Theorem \ref{th:main},
when $a \,\leq\, 0.2$ and $\lambda \,\leq\, 0.2 a$,
unless tolerance settings are forced to minimum values.
The true trajectory comes to much smaller population densities
and also spend more time at these very low population abundances.
Therefore, one has to be careful, since such results would give, e.g.,
a far to good picture of the populations chances to survive from any perturbation.

\begin{figure}[h]
\begin{center}
\includegraphics[scale = 0.9]{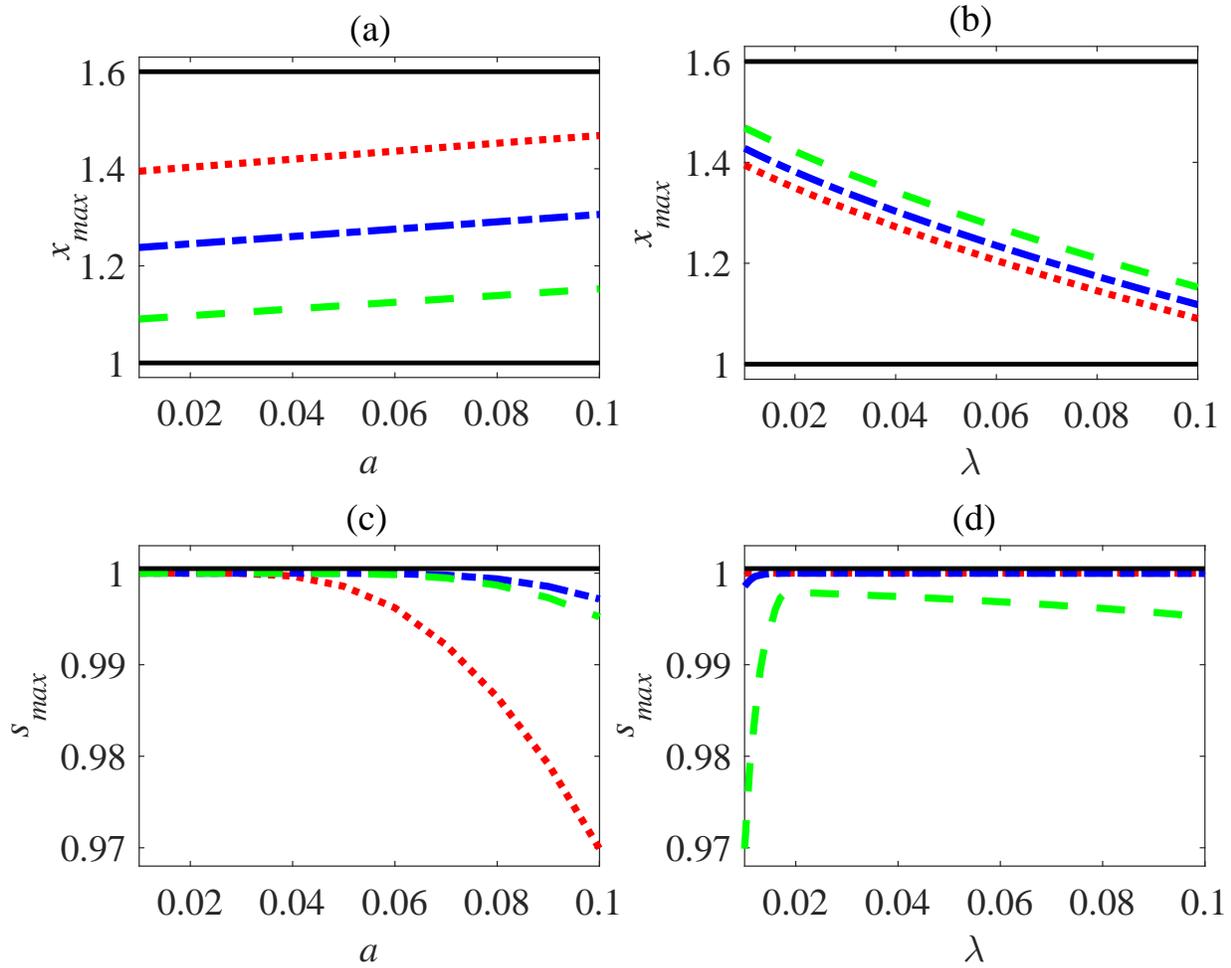}
\caption{The maximal $x$- and $s$-values as functions of the parameters $a$ and $\lambda$.
(a) and (c): $\lambda \,=\, 0.1$ (green, dashed), $\lambda \,=\, 0.05$ (blue, dashdot), $\lambda \,=\, 0.01$ (red, dotted).
(b) and (d): $a \,=\, 0.1$ (green, dashed), $a \,=\, 0.05$ (blue, dashdot), $a \,=\, 0.01$ (red, dotted).
Black solid lines show analytical estimates for $x_{max}$ and $s_{max}$ from Theorem \ref{th:main}.}
\label{fig:xmaxsmax}
\end{center}
\end{figure}

\begin{figure}[h]
\begin{center}
\includegraphics[scale = 0.85]{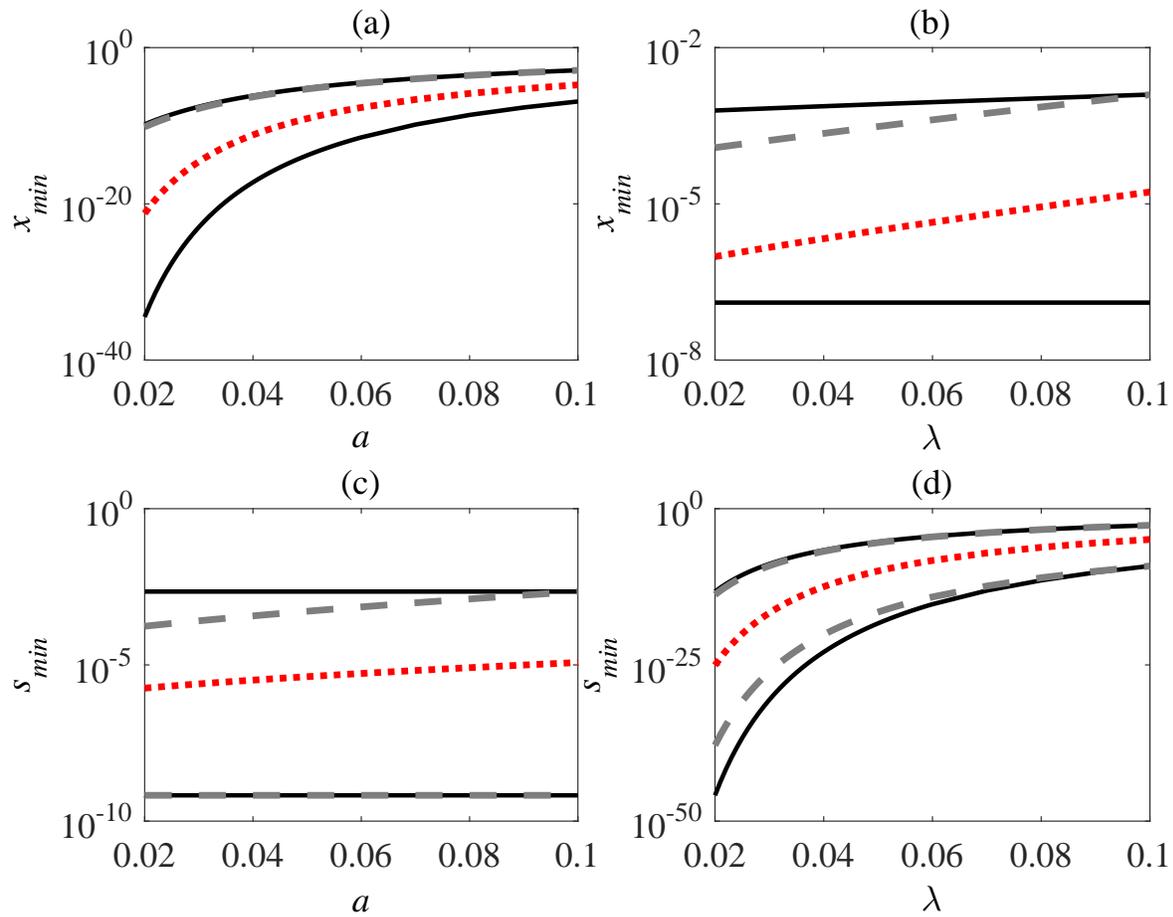}
\caption{The minimal $x$- and $s$-values as functions of the parameters $a$ and $\lambda$ (red, dotted):
(a) and (c): $\lambda \,=\, 0.1$.
(b) and (d): $a \,=\, 0.1$.
Grey dashed curves show analytical estimates for  $x_{min}$ and $s_{min}$ given in Theorem \ref{th:main} using
$\kappa_1$, $\kappa_2$ and $\kappa_3$,
while black solid curves show estimates produced by using the bounds of
$\kappa_1$, $\kappa_2$ and $\kappa_3$ given in Theorem \ref{th:main}.}
\label{fig:xminsmin}
\end{center}
\end{figure}

In Figure \ref{fig:xmaxsmax} the maximal $x$- and $s$-values, $x_{max}$ and $s_{max}$,
for the unique limit cycle are plotted as functions of the parameters $a$ and $\lambda$,
together with the analytical estimates given in Theorem \ref{th:main}.
Similarly, in Figure \ref{fig:xminsmin} the minimal $x$- and $s$-values, $x_{min}$ and $s_{min}$, for the unique limit cycle
are plotted.
The analytical estimates for $x_{min}$ and $s_{min}$ is produced by using the corresponding estimate for the maximal $x$-value, $x_{max} \,=\, 1.6$.
%

\begin{figure}[h]
\begin{center}
\includegraphics[height = 6cm, width = 7cm]{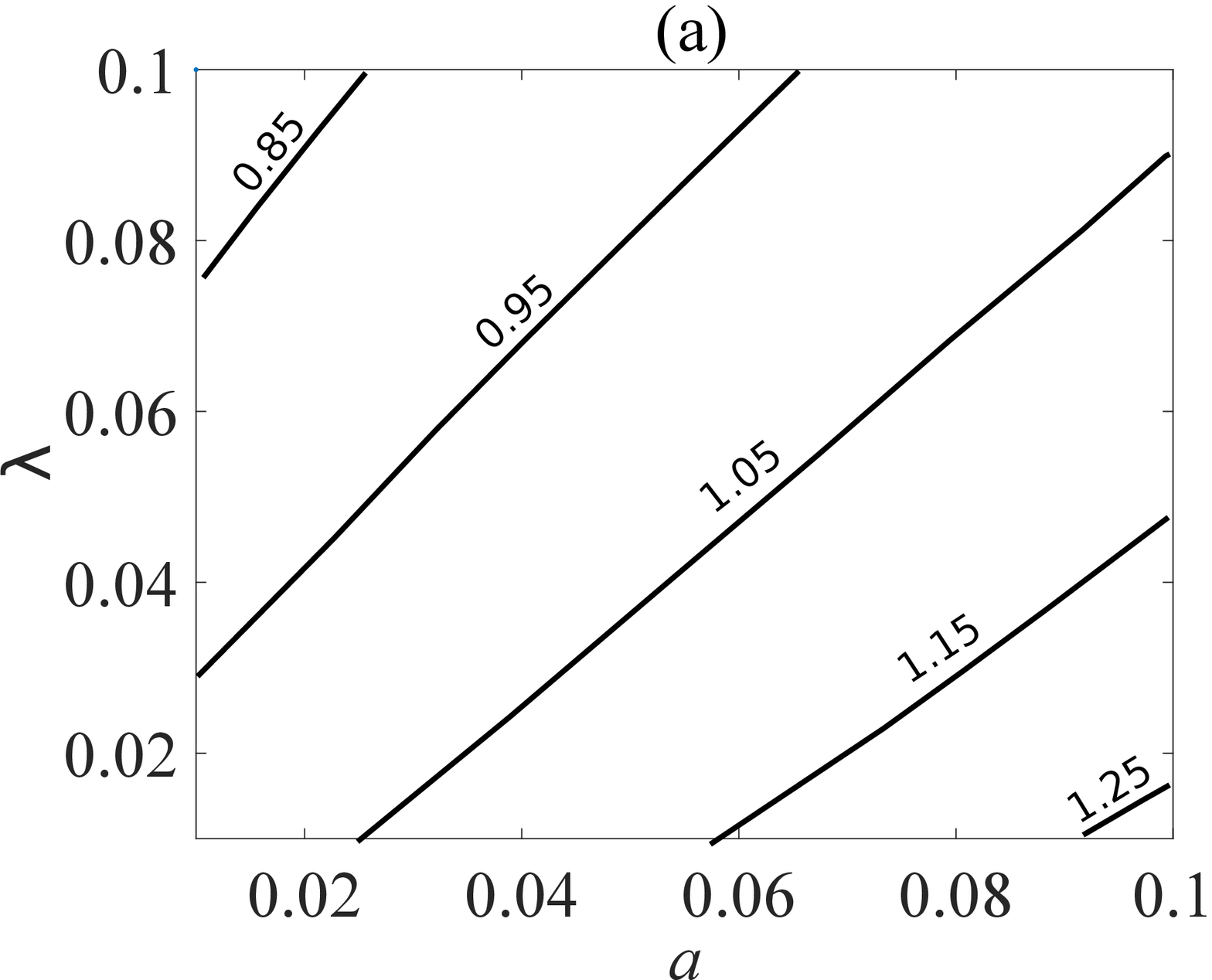}
\includegraphics[height = 6cm, width = 7cm]{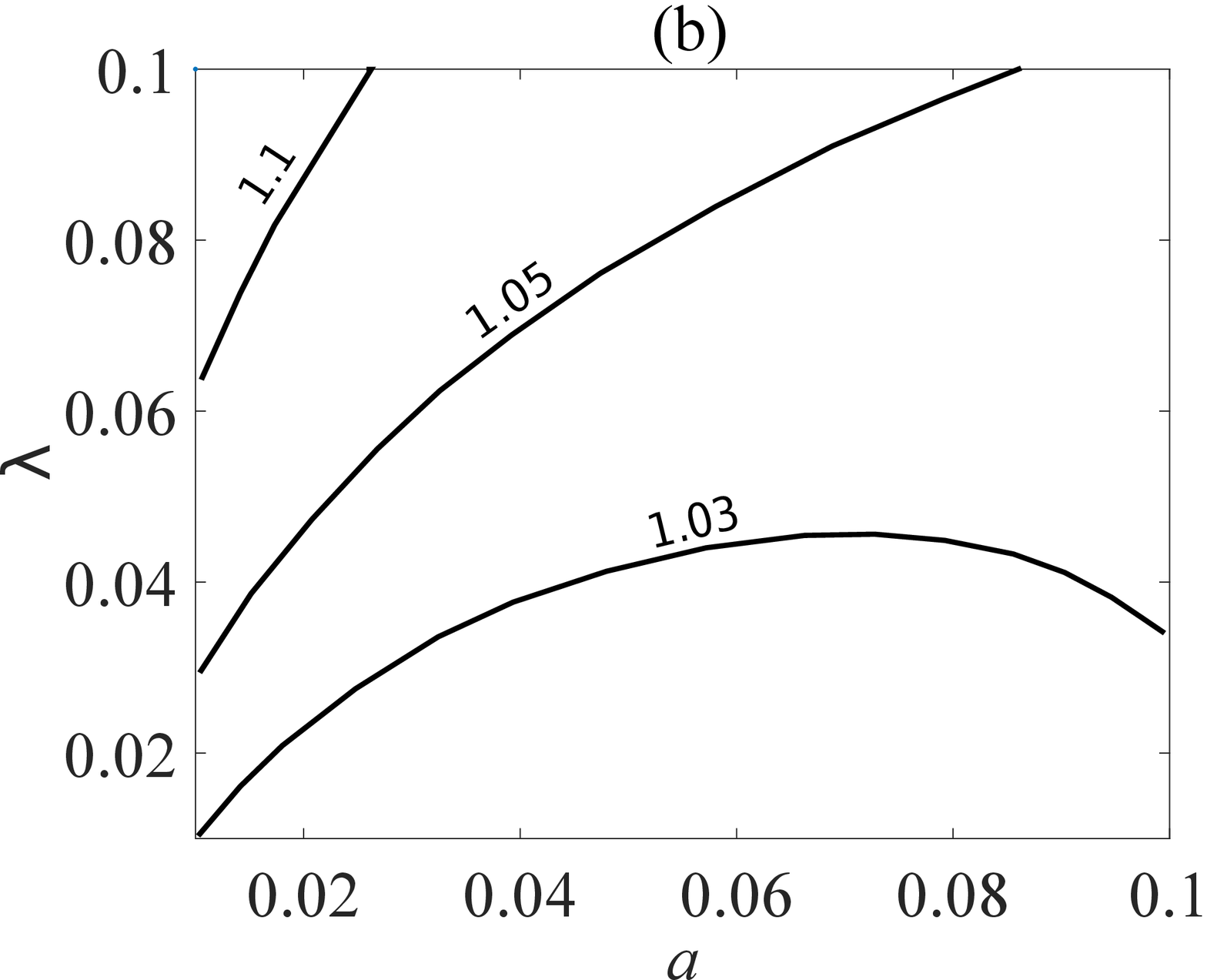}
\caption{Level curves of (a) $\tau_s$ and (b) $\tau_x$ as functions of $a$ and $\lambda$.
Observe the nonlinear steps between curves for $\tau_x$.
The function $\tau_s$ is close to one when $a \approx \lambda$,
while $\tau_x$ is close to one when $\lambda \approx 0.01$.}
\label{fig:levelcurves}
\end{center}
\end{figure}

Suppose now that $\left(x_{max},\lambda\right)$, $\left(h(s_{min}),s_{min}\right)$ and $\left(x_{min},\lambda\right)$
are points on the simulated limit cycle
and let $\tau_s$ and $\tau_x$ be such that
$$
s_{min} \,=\, \exp\left({-\frac{x_{max}}{\tau_s\lambda}}\right)
\qquad \text{and} \qquad
x_{min} \,=\, \exp\left({-\frac{x_{max}}{\tau_x a}}\right).
$$
We can thus say that $\tau_s$ and $\tau_x$ are measures of how good the approximations
$$
s_{min} \,\approx\, \exp\left({-\frac{x_{max}}{\lambda}}\right)
\qquad \text{and} \qquad
x_{min} \,\approx\, \exp\left({-\frac{x_{max}}{a}}\right),
$$
stated in Remark \ref{re:approx_main}, are.
Figure \ref{fig:levelcurves} shows level curves of the functions $\tau_s$ and $\tau_x$
in the $a \lambda$-plane for $a,\lambda \in (0.01,0.1)$.
From Figure \ref{fig:levelcurves} (a) we can observe that the approximation for $s_{min}$ is good for $a \approx \lambda$,
while Figure \ref{fig:levelcurves} (b) shows that the approximation for $x_{min}$ is good for $\lambda \approx 0.01$.

\begin{figure}[h]
\begin{center}
\includegraphics[scale = 0.85]{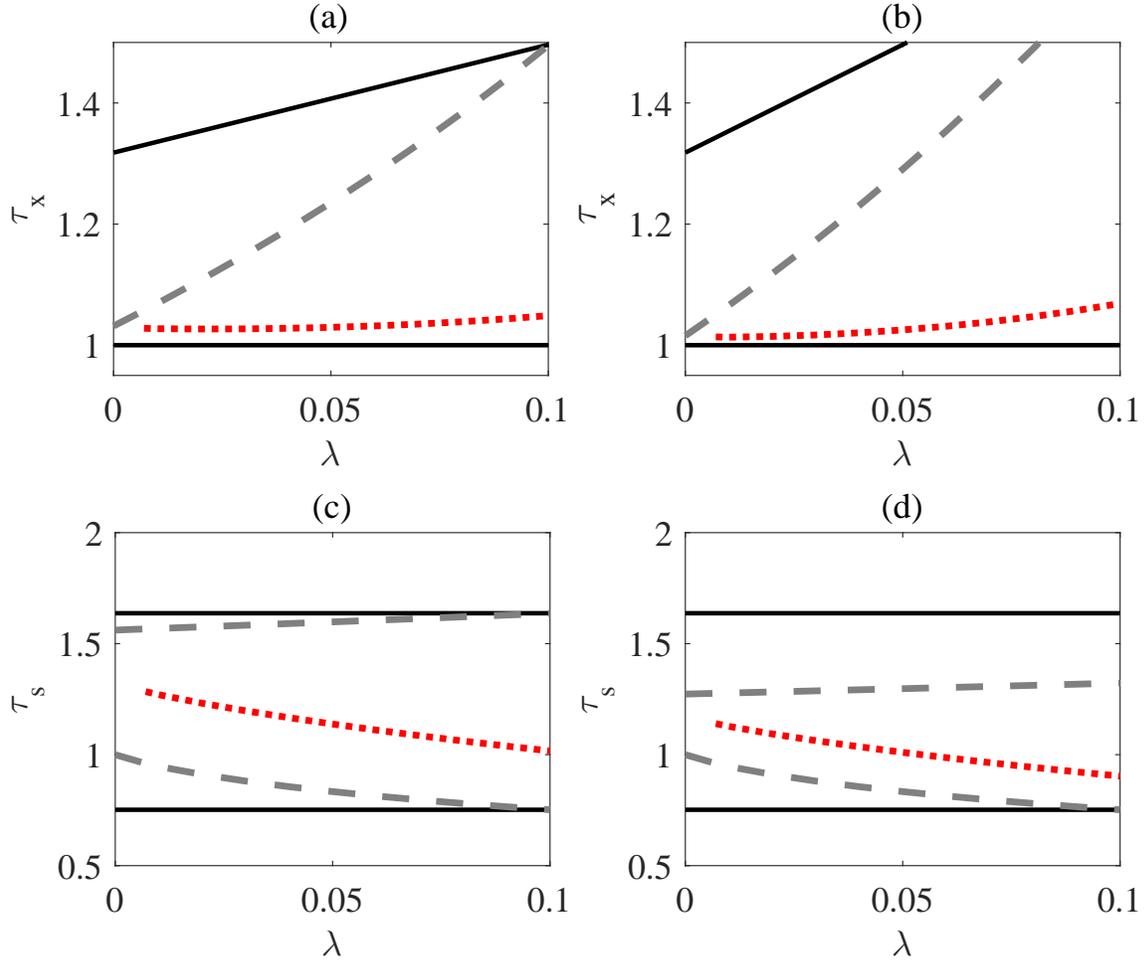}
\caption{The functions $\tau_s$ and $\tau_x$ as functions $\lambda$ (red, dotted):
(a) and (c) $a \,=\, 0.1$, (b) and (d) $a \,=\, 0.05$.
Grey dashed curves show analytical estimates for $\tau_x$ and $\tau_s$ given by
$\kappa_1$, $\kappa_2$ and $\kappa_3$ in Theorem \ref{th:main},
while the black solid curves show analytical estimates produced by using the bounds of
$\kappa_1$, $\kappa_2$ and $\kappa_3$ in Theorem \ref{th:main}.
}
\label{fig:xminsmin-tau}
\end{center}
\end{figure}

We end this section by plotting the functions $\tau_x$ and $\tau_s$,
for small values of $a$,
as functions of $\lambda$ in Figure \ref{fig:xminsmin-tau} together
with the analytical estimates for $\tau_x$ and $\tau_s$  given by
$\kappa_1$, $\kappa_2$ and $\kappa_3$ in Theorem \ref{th:main}.
As $\lambda$ and $\frac{\lambda}{a}$ approaches zero,
the lower estimate for $\tau_s$ approaches 1 ($\kappa_2 \to 1$)
while $\kappa_1$ and $\kappa_3$, giving upper estimates of $\tau_x$ and $\tau_s$,
stays a bit away from 1 for all $\lambda$.
%


.\newpage
.\newpage.\newpage.\newpage.\newpage


\bibliographystyle{amsalpha}

\end{document}